\documentclass[11pt, a4paper]{amsart}

\usepackage[T1]{fontenc}
\usepackage[utf8]{inputenc}
\usepackage[english]{babel}
\usepackage[leqno]{amsmath}
\usepackage{amssymb}
\usepackage{tabularx}
\usepackage{a4wide}
\usepackage{url}
\usepackage{tikz-cd}
\usepackage{csquotes}
\usepackage{enumitem}
\usepackage{mathtools}
\usepackage{manfnt}
\usepackage{mathrsfs}
\usepackage{aurical}

\usepackage[
  maxnames=4,
  maxalphanames=4,
  style=alphabetic,
  url = false,
  doi = false,
  isbn = false,
  ]{biblatex}
\renewbibmacro{in:}{}

\usepackage{todonotes} 

\definecolor{cadmiumgreen}{rgb}{0.0, 0.42, 0.24}
\definecolor{darkred}{rgb}{.85,0,0}
\usepackage[
  colorlinks,
  linkcolor=darkred,
  urlcolor=darkred,
  citecolor=cadmiumgreen,
  pdfstartview ={FitV},
]{hyperref}
\hypersetup{
  pdfauthor={Omid Amini, Matthieu Piquerez},
  pdftitle={Tropical Feichtner-Yuzvinsky and positivity criterion for fans},
  bookmarksdepth = 2}

\usepackage{scalerel}

\def\pathroundface#1#2[#3]{\path[#3] (I) -- (I#1) .. controls ($.3*(I#1)+.7*(I#1#1#2)$) and ($.3*(I#1#2)+.7*(I#1#1#2)$) .. (I#1#2) .. controls ($.3*(I#1#2)+.7*(I#2#1#2)$) and ($.3*(I#2)+.7*(I#2#1#2)$) .. (I#2) -- cycle}

\def\roundface#1#2[#3]{\pathroundface#1#2[fill,#3]}

\def\cornerface#1#2[#3]{\fill[#3] (I) -- (I#1) -- (I#1#1#2) -- (I#2#1#2) -- (I#2) -- cycle}

\usetikzlibrary{calc}
\tikzset{vcenter/.style={baseline={([yshift=-.8ex]current bounding box.center)}}}
\tikzset{dot/.style={insert path={node {\tikz[baseline=.6pt]\filldraw[black] (0,0) circle (1.2pt);}}}}

\setlist[itemize,1]{itemsep=\smallskipamount}
\setlist[enumerate,1]{itemsep=\smallskipamount, label=\textnormal{(\arabic*)}}

\newtheorem{thm}{Theorem}[section]
\newtheorem{lemma}[thm]{Lemma}
\newtheorem{prop}[thm]{Proposition}

\theoremstyle{definition}

\newtheorem{remm}[thm]{Remark}
\newenvironment{remark}
  {\pushQED{\qed}\remm}
  {\popQED\endremm}
\newtheorem{exx}[thm]{Example}
\newenvironment{example}
  {\pushQED{\qed}\exx}
  {\popQED\endexx}

\numberwithin{equation}{section}

\newcommand{\ie}{i.e.}

\newcommand{\resp}{resp.\ }

\renewcommand{\~}{\widetilde}
\renewcommand{\hat}{\widehat}
\newcommand{\Q}{\mathbb{Q}}
\newcommand{\Z}{\mathbb{Z}}
\newcommand{\N}{\mathbb{N}}
\newcommand{\R}{\ss{\mathbb{R}}}

\newcommand{\vsim}{\rotatebox{90}{$\sim$}}

\let\oldforall\forall
\renewcommand{\forall}{\oldforall\:}
\let\oldbigwedge\bigwedge
\renewcommand{\bigwedge}{{\textstyle\oldbigwedge\!}}

\renewcommand{\geq}{\geqslant}
\renewcommand{\leq}{\leqslant}
\renewcommand{\setminus}{\smallsetminus}

\makeatletter
\let\oldsum\sum
\renewcommand{\sum}{\@ifnextchar_\@mysum\oldsum}
\def\@mysum_#1{\oldsum_{\substack{#1}}}
\let\oldbigoplus\bigoplus
\renewcommand{\bigoplus}{\@ifnextchar_\@mybigoplus\oldbigoplus}
\def\@mybigoplus_#1{\oldbigoplus_{\substack{#1}}}
\let\oldprod\prod
\renewcommand{\prod}{\@ifnextchar_\@myprod\oldprod}
\def\@myprod_#1{\oldprod_{\substack{#1}}}
\makeatother

\makeatletter
\let\oldnu\nu
\newlength{\heightnu}
\newlength{\depthnu}
\def\nu#1_#2{{\settoheight{\heightnu}{\hbox{$#2$}}\settodepth{\depthnu}{\hbox{$#2$}}\oldnu\rule[\depthnu-3pt]{0pt}{1pt}#1_{\!#2}}}
\makeatother

\makeatletter
\let\@oldinfty\infty
\newcommand{\@sminfty}{{\scaleto{\@oldinfty}{2.8pt}}} 
\renewcommand{\infty}{{\mathchoice%
  {\displaystyle{\@oldinfty}}%
  {\textstyle{\@oldinfty}}%
  {\hspace{-1pt}\scriptstyle{\@sminfty}}%
  {\hspace{-1pt}\scriptscriptstyle{\@sminfty}}}
}
\makeatother

\newcommand{\rquot}[2]{#1\big/#2}
\newcommand{\rest}[1]{\raisebox{-1pt}{$\vert$}_{#1}}
\newcommand{\simto}{\xrightarrow{\raisebox{-3pt}[0pt][0pt]{\small$\hspace{-1pt}\sim$}}}
\newcommand{\longsimto}{\xrightarrow{\ \raisebox{-3pt}[0pt][0pt]{\small$\hspace{-1pt}\sim$\ }}}
\newcommand{\id}{\mathrm{id}} 
\newcommand{\dual}{\star}
\newcommand{\bul}{\bullet}

\let\hom\relax
\DeclareMathOperator{\hom}{Hom} 

\DeclareMathOperator{\coker}{coker}

\newcommand{\gr}{\mathrm{gr}} 
\renewcommand{\i}{\mathrm i} 
\renewcommand{\d}{\mathrm d} 
\newcommand{\BM}{{^{\scaleto{\mathrm{BM}}{3.6pt}}}}
\newcommand{\cubC}{C_{\ss{}_\square}} 
\newcommand{\hcubC}{C^{^\square}} 
\newcommand{\cubd}{\d_{\ss{}_\square}} 
\DeclareMathOperator{\PD}{PD} 
\DeclareMathOperator{\sign}{sign} 

\newcommand{\T}{\ss{\mathbb T}} 
\renewcommand{\P}{\ss{\mathbb P}} 
\newcommand{\TP}{{\T\P}} 

\newcommand{\SF}{\mathbf F} 
\newcommand{\RMod}{{\R_+}} 
\newcommand{\opencone}[1]{\ss{\mathring{#1}}} 

\newcommand{\conezero}{{\underline0}}
\newcommand{\e}{\mathfrak e} 

\newcommand{\nvect}{\ss{\mathfrak n}} 
\newcommand{\chart}{\~} 
\newcommand{\cube}{\ssbis[0pt][.5pt]{\scaleto{\square}{6pt}}!} 
\DeclareMathOperator{\sed}{sed} 

\newcommand{\shiftcomp}[2][0]{{}\mkern#1mu\overline{\mkern-#1mu#2}}
\newcommand{\comp}[1]{\if#1X \shiftcomp[3]{#1}\else\if#1Z \shiftcomp[3]{#1} \else \shiftcomp{#1}\fi\fi} 

\newcommand{\suppaux}[2]{\scalebox{1}[1.4]{$#1\lvert$}#2\scalebox{1}[1.4]{$#1\rvert$}}
  \newcommand{\supp}[1]{\mathpalette\suppaux{#1}}
\newcommand{\dimsaux}[2]{\raisebox{.2ex}{\scalebox{1}[.8]{$#1\lvert$}}#2\raisebox{.2ex}{\scalebox{1}[.8]{$#1\rvert$}}}
  \newcommand{\dims}[1]{\mathpalette\dimsaux{#1}}
\newcommand{\subface}{\prec}
\newcommand{\ssface}{\mathbin{\mathchoice
  {\subface\!\!\!\cdot}%
  {\subface\!\!\!\cdot}%
  {\subface\!\cdot}%
  {\subface\!\cdot}%
}} 
\newcommand{\supface}{\succ}
\newcommand{\ssupface}{\mathbin{\mathchoice
  {\cdot\!\!\!\supface}%
  {\cdot\!\!\!\supface}%
  {\cdot\!\supface}%
  {\cdot\!\supface}%
}}

\newcommand{\MW}{\mathrm{MW}} 
\newcommand{\W}{\mathrm{W}} 
\newcommand{\x}{\ss{\textsc{x}}}

\newcommand{\Ma}{{\scalebox{1.12}{$\mathfrak m$}}}
\DeclareMathOperator{\cl}{cl} 

\newcommand{\E}{\textnormal{\textsf{E}}}
\newcommand{\Epnop}{\E}
\newcommand{\Ep}[1]{\prescript{}{#1}\Epnop} 
\newcommand{\EpI}[1]{\prescript\downarrow{#1}\Epnop} 
\newcommand{\EpII}[1]{\prescript\rightarrow{#1}\Epnop} 
\newcommand{\F}{\textnormal{\textsf{F}}}
\newcommand{\Fpq}[2]{\prescript{#1}{#2}\F} 
\DeclareMathOperator{\Tot}{Tot} 

\newcommand{\Ext}[2]{\mathrm{Ext}^1\bigl(#1, #2\bigr)} 
\newcommand{\contract}{\kappa}

\RenewDocumentCommand{\ss}{O{0pt} O{0pt} O{.8} m e{_^}}{
  #4%
  \IfValueT{#5}{
    \sb{\hspace{#1}\scaleobj{#3}{#5}}
  }
  \IfValueT{#6}{
    \sp{\hspace{#2}\scaleobj{#3}{#6}}
  }
}

\NewDocumentCommand{\ssbis}{O{0pt} O{0pt} O{.8} m t! e{_^}}{
  #4%
  \IfValueT{#6}{
    \IfBooleanTF{#5}{\sb{\hspace{#1}\scaleobj{#3}{#6}}}{\sb{#6}}
  }
  \IfValueT{#7}{
  \IfBooleanTF{#5}{\sp{\hspace{#2}\scaleobj{#3}{#7}}}{\sp{#7}}
}
}

\newcommand{\ssomega}{\ss\omega}
\newcommand{\sspi}{\ss\pi}

\newcommand{\sse}{\ss\e}

\newcommand{\ssTP}{\ss\TP}
\newcommand{\ssi}{\ss\i}
\newcommand{\ssrho}{\ss{\rho}}
\newcommand{\sscompsigma}{\ss[-1pt]{\comp\sigma}}
\newcommand{\sssigma}{\ss[-1pt]{\sigma}}
\newcommand{\sstau}{\ss{\tau}}
\newcommand{\sseta}{\ss{\eta}}

\newcommand{\ssM}{\ss{M}}
\newcommand{\ssN}{\ss{N}}
\newcommand{\ssH}{\ss{H}}
\newcommand{\ssvarpi}{\ss{\varpi}}
\newcommand{\ssnu}{\ss[-1pt]{\oldnu}}

\newcommand{\ssSigma}{\ss{\Sigma}}
\newcommand{\compSigma}{\ss{\comp\Sigma}}

\newcommand{\ssinfty}{\ss{\infty}}
\newcommand{\compcone}[1]{\ssbis[-2pt][.5pt]{\comp{#1}}!}
\newcommand{\ssW}{W}
\newcommand{\ssalpha}{\ss{\alpha}}
\newcommand{\ssx}{\ss{x}}
\newcommand{\ssf}{\ss{f}}

\newcommand{\ssa}{\ss{a}}
\newcommand{\ssb}{\ss{b}}
\newcommand{\ssc}{\ss{c}}
\newcommand{\sshata}{\ss{\hat a}}

\makeatletter
\@namedef{subjclassname@2020}{\textup{2020} Mathematics Subject Classification}
\makeatother

\newcommand{\msc}[1]{\href{http://www.ams.org/msc/msc2020.html?t=&s=#1}{#1}}

\begin{document}

\allowdisplaybreaks

\title{Tropical Feichtner-Yuzvinsky and positivity criterion for fans}

\author{Omid Amini}
\address{CNRS - CMLS, \'Ecole polytechnique, Institut polytechnique de Paris.}
\email{\href{omid.amini@polytechnique.edu}{omid.amini@polytechnique.edu}}

\author{Matthieu Piquerez}
\address{LS2N, Inria, Nantes Université}
\email{\href{matthieu.piquerez@univ-nantes.fr}{matthieu.piquerez@univ-nantes.fr}}


\date{\today}

\begin{abstract} We prove that the Chow ring of any simplicial fan is isomorphic to the middle degree part of the tropical cohomology ring of its canonical compactification. Using this result, we prove a tropical analogue of Kleiman's criterion of ampleness for fans.

In the case of tropical fans that are homology manifolds, we obtain an isomorphism between the Chow ring of the fan and the entire tropical cohomology of the canonical compactification. When applied to matroids, this provides a new representation of the Chow ring of a matroid as the cohomology ring of a projective tropical manifold.
\end{abstract}

\maketitle

\setcounter{tocdepth}{1}

\tableofcontents

\section{Introduction}

Feichtner and Yuzvinsky provide in~\cite{FY04} a description of the Chow ring of a wonderful compactification of the complement of an arrangement of hyperplanes in terms of the Chow ring of the toric variety that underlines the compactification. Over the field of complex numbers, combined with the work of De Concini and Procesi~\cite{DP95-2}, this leads to a combinatorial description of the cohomology of the wonderful compactification. Our aim in this paper is to prove a very general form of these results in the framework of tropical geometry. This plays a key role in our companion work which develops a Hodge theory for Kähler tropical varieties.

\subsection{Overview of the main results}

Let $N$ be a free abelian group of finite rank and denote by $\ssN_\R$ the vector space generated by $N$. Let $\Sigma$ be a simplicial fan in $\ssN_{\R}$ rational with respect to $N$, and denote by $\comp\Sigma$ its canonical compactification obtained by taking the closure of $\Sigma$ in the tropical toric variety $\ssTP_\Sigma$. Denote by $A^\bul(\Sigma, \Z)$ and $A^\bul(\Sigma, \Q)$ the Chow ring of $\Sigma$ with integral and rational coefficients, and by $H^{\bul,\bul}(\comp\Sigma, \Z)$ and $H^{\bul,\bul}(\comp\Sigma, \Q)$ the tropical cohomology groups of $\comp\Sigma$ with integral and rational coefficients, respectively. We refer to Section~\ref{sec:prel} for a reminder of these definitions.

\smallskip
Our first theorem states as follows.

\begin{thm}[Tropical Feichtner-Yuzvinsky for fans] \label{thm:ring_morphism}
For any integer $p$, there is an isomorphism
\[ \begin{array}{ccc}
\Psi\colon H^{p,p}(\comp\Sigma, \Q) & \longsimto & A^p(\Sigma, \Q).
\end{array} \]
They induce together a ring morphism $H^{\bul,\bul}(\comp\Sigma, \Q) \to A^\bul(\Sigma, \Q)$ by mapping $H^{p,q}(\comp\Sigma, \Q)$ to zero in the bidegree $p\neq q$. Moreover, $H^{p,q}(\comp\Sigma, \Q)$ is trivial for $p < q$ and for $p>q=0$.

If\/ $\Sigma$ is unimodular and saturated, these statements hold with $\Z$-coefficients.
\end{thm}

We refer to Section~\ref{subsec:fans} for the definition of saturation property. The proof given in Section~\ref{sec:first_main_theorem} provides an explicit description of the application $\Psi$. In Sections~\ref{sec:ring_morphism_non_saturated} and~\ref{sec:non_unimodular_projective_fan} we will provide examples which show that both the saturation and unimodularity assumptions are in general needed when dealing with integral cohomology.

The above result has its source of motivation in the recent development of combinatorial Hodge theory~\cite{AHK, Ard22, Baker18, Huh22, Okou22}. To any matroid $\Ma$ is associated a unimodular fan $\ssSigma_\Ma$ called the Bergman fan of the matroid~\cite{AK06}. The Chow ring of $\Ma$ is by definition the Chow ring of the corresponding Bergman fan. It behaves as the cohomology ring of a smooth projective complex variety, with remarkable properties, although, in general, when the matroid is non-realizable, there is no projective variety associated to $\Ma$.

From Theorem~\ref{thm:ring_morphism}, we deduce the following identification of the Chow ring of a matroid as the cohomology of a projective tropical manifold.

\begin{thm}\label{thm:Hodge_isomorphism_matroids}
Let $\ssSigma_\Ma$ be the Bergman fan of a matroid $\Ma$ and denote by $\ss{\comp\Sigma}_\Ma$ its canonical compactification. Denote by $A^\bul(\Ma, \Z)$ the Chow ring of\/ $\Ma$ with integral coefficients. We have an isomorphism of rings $A^\bul(\Ma, \Z) \simto H^{\bul, \bul}(\ss{\comp\Sigma}_\Ma, \Z)$.
\end{thm}

We will deduce the above theorem as a special case of the following more general result. A tropical orientation of a rational fan $\Sigma$ of pure dimension $d$ is an integer valued map
\[ \omega\colon \ssSigma_d \to \Z\setminus\{0\} \]
that verifies the so-called balancing condition. Namely, for any cone $\tau$ in $\Sigma$ of codimension one, denoting by $\ssN_\tau$ the lattice of full rank in the vector subspace of $\ssN_\R$ generated by $\tau$, we have the vanishing of the sum
\[ \sum_{\sigma \supset \tau} \omega(\sigma)\sse_{\sigma}^{\tau} =0 \]
in the quotient lattice $\rquot{N}{N_\tau}$. Here, the sum runs over facets $\sigma$ of $\Sigma$ that contain $\tau$, and $\e_{\sigma}^{\tau}$ is the generator of the quotient $\rquot{(\sigma \cap N)}{(\tau \cap N)} \simeq \Z_{\geq 0}$.

A tropical fan in $N_{\R}$ is a pair $(\Sigma, \ssomega_\Sigma)$ consisting of a pure dimensional rational fan $\Sigma$ and a tropical orientation $\ssomega_\Sigma$ as above. The tropical orientation leads to the definition of a fundamental class. We say that a tropical fan is a tropical homology manifold if the tropical cohomology with rational coefficients of any open subset of $\supp{\Sigma}$ with induced orientation verifies Poincaré duality. This has been thoroughly studied in~\cite{Aks21,AP-homology}. If the same holds with integral coefficients, then we precise that the tropical fan is a tropical homology manifold with integral coefficients.

The Bergman fan of a matroid $\Ma$ is pure dimensional and has a natural tropical orientation that gives value one to each facet. It is shown in~\cite{JSS19} that any open set in the Bergman fan of a matroid verifies Poincaré duality, that is, $\ssSigma_\Ma$ is a tropical homology manifold.

Theorem~\ref{thm:Hodge_isomorphism_matroids} is a special case of the following more general result.

\begin{thm}\label{thm:Hodge_isomorphism}
Let $\Sigma$ be a simplicial tropical fan in $\ssN_{\R}$. Suppose that in addition, $\Sigma$ is a tropical homology manifold. Then, we get an isomorphism of rings $A^\bul(\Sigma, \Q) \simto H^{\bul, \bul}(\comp\Sigma, \Q)$. In particular, $A^\bul(\Sigma, \Q)$ verifies Poincaré duality.

If\/ $\Sigma$ is unimodular, saturated, and a tropical homology manifold with integral coefficients, the statement holds over $\Z$. In this case, $A^\bul(\Sigma, \Z)$ is torsion-free and verifies Poincaré duality.
\end{thm}

Combined with~\cite{AP-hodge-fan}, these results show that the recently developed Hodge theory for tropical fans concern the cohomology of special projective tropical varieties, those of the form~$\comp\Sigma$. In our companion work~\cite{AP-tht}, these results are used to introduce Kähler tropical varieties and develop a Hodge theory for them.

\smallskip
The following dual version of Theorem~\ref{thm:ring_morphism} does not use the saturation hypothesis. Denote by $\MW_p(\Sigma, \Z)$ the group of integral valued Minkowski weights of dimension $p$ and consider the cycle class map $\cl\colon \MW_p(\Sigma, \Z) \to H_{p,p}(\comp\Sigma, \Z)$.
We refer to Section~\ref{sec:MW} for the definition.

\begin{thm} \label{thm:Hodge_isomorphism_dual}
Let $\Sigma$ be a unimodular fan in $\ssN_{\R}$. The cycle class map $\MW_p(\Sigma, \Z) \to H_{p,p}(\comp\Sigma, \Z)$ is an isomorphism. Moreover, $H_{p,q}(\comp\Sigma, \Z)$ is trivial in bidegree $(p,q)$ in the case $p < q$ and $p > q= 0$.
\end{thm}

We note that for rational simplicial fans, the statement of the above theorem with rational coefficients is a direct consequence of Theorem~\ref{thm:Hodge_isomorphism} and of the duality between Minkowski weights and Chow groups, Theorem~\ref{thm:duality_chow_mw}.

We get the following corollary.

\begin{thm}[Hodge conjecture for compactifications of tropical fans] \label{thm:Hodge_conjecture}
Let\/ $\Sigma$ be a unimodular tropical fan which is a tropical homology manifold with integral coefficients. Then, to each element $\alpha$ in $H^{p,p}(\comp\Sigma, \Z)$ we can associate a tropical cycle of codimension $p$ whose class in $H_{d-p,d-p}(\comp\Sigma, \Z)$ is the Poincaré dual of $\alpha$. Moreover, for $p \neq q$, $H_{p,q}(\comp\Sigma, \Z)$ is trivial.

The same statements are true with rational coefficients for a simplicial tropical fan $\Sigma$ which is a tropical homology manifold.
\end{thm}

Using the above results, we formulate a tropical analogue of Kleiman's criterion of ampleness~\cite{Kle66}. An element $\alpha$ in $A^1(\Sigma, \Q)$ is called ample if it is associated to a conewise linear function on $\Sigma$ that is strictly convex around each cone of $\Sigma$. We say that a tropical variety is effective if its underlying tropical orientation gets positive values. We prove the following numerical characterization of ampleness for fans in Section~\ref{sec:kleiman}.

\begin{thm}[Numerical criterion of ampleness for fans]\label{thm:positivity}
Let $\Sigma$ be a rational simplicial fan. A class $\alpha$ in $A^1(\Sigma,\Q)$ is ample if and only if, viewed in $H^{1,1}(\comp\Sigma, \Q)$, has positive pairing with the class in $H_{1,1}(\comp\Sigma, \Q)$ of any effective tropical curve in $\comp\Sigma$.
\end{thm}

\subsection{An alternative description of the cohomology}

The canonical compactification $\comp\Sigma$ admits a natural decomposition into hypercubes, see Figure~\ref{fig:canonical_compactification}. In order to prove Theorem~\ref{thm:ring_morphism}, we use this decomposition and establish a new way of computing the cohomology of $\comp\Sigma$ that we hope might be of independent interest. We briefly discuss this here.

For each cone $\sigma \in \Sigma$, the canonical compactification $\comp\sigma$ of $\sigma$ is a subset of $\comp\Sigma$. It has the form of a hypercube of dimension equal to that of $\sigma$. Denote by $\ssinfty_\sigma$ the point in $\comp\sigma$ that is diagonally opposite to the origin. For non-negative integers $p,q$, define
\begin{equation}\label{eq:definition_cubical_complex}
\cubC^{p,q}(\comp\Sigma, \Z) \coloneqq \bigoplus_{\dims\sigma = q} \SF^{p-q}(\ssinfty_\sigma, \Z)
\end{equation}
where $\SF^{k}(\ssinfty_\sigma, \Z)$ is the coefficient group used in the definition of tropical cohomology. It coincides with tropical cohomology in degree $k$ of the star fan $\ssSigma^\sigma$, see Section~\ref{sec:prel} for the definition.

For $\tau$ a face of codimension one in $\sigma \in \Sigma$, we define a natural map
\[ \SF^{k}(\ssinfty_\tau, \Z)  \longrightarrow  \SF^{k-1}(\ssinfty_\sigma, \Z).\]
Viewing these as boundary maps and combining them together, we obtain for each integer $p$ a cochain complex
\[ \cubC^{p,\bullet}(\comp\Sigma, \Z)\colon \quad \dots\longrightarrow \cubC^{p, q-1}(\comp\Sigma, \Z) \xrightarrow{\ \cubd^{q-1}\ } \cubC^{p,q}(\comp\Sigma, \Z) \xrightarrow{\ \cubd^{q}\ } \cubC^{p,q+1}(\comp\Sigma, \Z) \longrightarrow \cdots \]

\begin{thm} \label{thm:cohomology_compactification}
Let $\Sigma$ be a unimodular fan. The cohomology of\/ $\bigl(\cubC^{p,\bul}(\comp\Sigma, \Z), \cubd \bigr)$ with integral coefficients is isomorphic to $H^{p,\bul}(\comp\Sigma, \Z)$. The same statement holds true with rational coefficients for simplicial rational fans.
\end{thm}

The proof of the above theorem is based on the homological properties of the \emph{fine double complex} $\Ep{p}^{\bul,\bul}$ introduced in Section~\ref{sec:fine_double_complex}.

A deeper exploitation of the fine double complex leads to the proof of the following characterization of tropical homology manifolds. We state it for tropical fans.

\begin{thm}[Alternate characterization of tropical homology manifolds] \label{thm:smoothness_alternate}
  A unimodular tropical fan $\Sigma$ is a tropical homology manifold with integral coefficients if and only if, for any face $\sigma\in\Sigma$, the cohomology with integral coefficients of\/ $\ss{\comp\Sigma}^\sigma$ verifies Poincaré duality. Equivalently, if and only if, for each $\sigma \in \Sigma$,
  \begin{itemize}
  \item the Chow ring $A^\bul(\ssSigma^\sigma, \Z)$ of the star fan $\ssSigma^\sigma$ verifies Poincaré duality, and
  \item all the cohomology groups $H^{p,q}(\ss{\comp{\Sigma}}^\sigma, \Z)$ for $p>q$ are zero.
\end{itemize}

The same statements hold true for simplicial tropical fans $\Sigma$ with rational coefficients.
\end{thm}

\subsection{Generalization to non-rational simplicial fans}

Although we assume for the sake of simplicity and clarity of the exposition that the fans are rational, it is possible to extend the set-up to non-rational simplicial fans. We discuss the required adjustments in Section~\ref{sec:nonrational}.

\subsection{Applications}

The results of this paper play an important role in our companion work. We briefly discuss them here.

Tropical Hodge theory in the global setting is the subject of our work~\cite{AP-tht}. Theorems~\ref{thm:ring_morphism} and~\ref{thm:Hodge_isomorphism}, combined with Steenbrink-Tropical comparison theorem and tropical Deligne weight exact sequence, enables us to use Chow rings of tropical fans in the study of the cohomology of tropical varieties. Using this approach, we introduce Kähler tropical varieties and establish a Hodge theory for them. Theorem~\ref{thm:Hodge_isomorphism_matroids} connects combinatorial and tropical Hodge theories.

Theorem~\ref{thm:Hodge_isomorphism} plays as well a crucial role in our joint work with Aksnes and Shaw~\cite{AAPS}. In that paper, we provide a characterization of varieties for which tropicalization remembers the cohomology. This is motivated in part by the work of Deligne~\cite{Del-md}, in which he gives a Hodge-theoretic characterization of maximal degenerations of complex algebraic varieties, and by the recent work by Yang Li~\cite{Li20}, that connects SYZ conjecture in mirror symmetry to tropical geometry. In this regard, we develop in~\cite{AP23-MA} a differential calculus on tropical varieties that combine Chow rings of local tropical fans with real differential forms on the variety, formulate a tropical analogue of the Monge-Ampère equation, and study its solutions.

\subsection{Organization}

Section~\ref{sec:prel} contains preliminary definitions and and results. In Section~\ref{sec:fine_double_complex}, we introduce the fine double complex and prove Theorem~\ref{thm:cohomology_compactification}. In Section~\ref{sec:first_main_theorem}, we prove Theorems~\ref{thm:ring_morphism},~\ref{thm:Hodge_isomorphism_matroids},~\ref{thm:Hodge_isomorphism},~\ref{thm:Hodge_isomorphism_dual} and~\ref{thm:Hodge_conjecture}. The assertion that the morphism $\Psi$ in Theorem~\ref{thm:ring_morphism} is a morphism of rings is more subtle and is proved in Section~\ref{subsec:Ap_to_Hpp}. Theorem ~\ref{thm:smoothness_alternate} is proved in Section~\ref{sec:proof-smoothness-criterion}, and Kleiman's criterion is established in Section~\ref{sec:kleiman}. Concluding Section~\ref{sec:discussion} contains complementary results and examples related to the content of the paper.

\subsection*{Convention}

In the following, whenever there is a choice, we will present the proofs of the results stated in the theorems for the case of integer coefficients, assuming that the required conditions, if any, are satisfied. The proofs with rational coefficients (without those conditions) are similar. Moreover, when the coefficient is not explicitly given, it means that we work with integral coefficients.

\subsection*{Acknowledgments}

We thank Pierre-Louis Blayac for his help in proving Theorem~\ref{thm:positivity}.

\section{Preliminaries}\label{sec:prel}

In this section we collect basic notations and definitions that are used in the paper.

Throughout, $N$ will be a free $\Z$-module of finite rank and $M=N^\dual = \hom(N, \Z)$ will be the dual of $N$. We denote by $\ssN_{\R}$ and $\ssM_{\R} = \ssN_{\R}^\dual$ the corresponding real vector spaces.

All the cones appearing in this paper are strictly convex, \ie, they do not contain any line. For a rational polyhedral cone $\sigma$ in $\ssN_{\R}$, we use the notation $\ssN_{\sigma, \R}$ to denote the real vector subspace of $\ssN_{\R}$ generated by elements of $\sigma$ and set $\ssN_{\R}^\sigma \coloneqq \rquot{\ssN_{\R}}{\ssN_{\sigma, \R}}$. Since $\sigma$ is rational, we get natural lattices of full rank $\ssN_\sigma \subset \ssN_{\sigma, \R}$ and $\ssN^\sigma \subset \ssN_{\R}^\sigma$.

For the ease of reading, we adopt the following convention. We use $\sigma$ (or any other face of~$\Sigma$) as a superscript where referring to the quotient of some space by $\ssN_{\sigma, \R}$ or to the elements related to this quotient. In contrast, we use $\sigma$ as a subscript for subspaces of $\ssN_{\sigma,\R}$ or for elements associated to these subspaces.

If $\tau\subseteq \sigma$ are faces of $\Sigma$, we denote by $\sspi_{\tau \subface \sigma}$ both the projection maps $\ssN^\tau\to \ssN^\sigma$ and $\ssN^\tau_{\R}\to \ssN^\sigma_{\R}$.

We denote by $\T \coloneqq \R \cup \{\ssinfty\}$ the extended real line with the topology induced by that of $\R$ and a basis of open neighborhoods of infinity given by intervals $(a, \ssinfty]$ for $a\in \R$. Extending the addition of $\R$ to $\T$ by setting $\ssinfty + a = \ssinfty$ for all $a \in \T$, endows $\T$ with the structure of a topological monoid called the monoid of tropical numbers. We denote by $\T_+ \coloneqq \R_+ \cup\{\ssinfty\}$ the submonoid of non-negative tropical numbers with the induced topology. Both monoids admit a natural scalar multiplication by non-negative real numbers (setting $0\cdot\ssinfty=0$). Moreover, the multiplication map is continuous. As such, $\T$ and $\T_+$ can be seen as modules over the semiring $\R_+$. A polyhedral cone $\sigma$ in $\ssN_{\R}$ is another example of a module over $\R_+$.

For any natural number $p$, we set $[p] \coloneqq \{1, \dots, p\}$.

\subsection{Fans} \label{subsec:fans}

Let $\Sigma$ be a fan of dimension $d$ in $\ssN_{\R}$. The dimension of a cone $\sigma$ in $\Sigma$ is denoted by $\dims\sigma$. The set of $k$-dimensional cones of $\Sigma$ is denoted by $\ssSigma_k$, and elements of $\ssSigma_1$ are called \emph{rays}. We denote by $\conezero$ the cone $\{0\}$. Any $k$-dimensional cone $\sigma$ in $\Sigma$ is determined by its set of rays in $\ssSigma_1$. The \emph{support} of $\Sigma$ denoted $\supp \Sigma$ is the closed subset of $\ssN_{\R}$ obtained by taking the union of the cones in $\Sigma$, we call it a \emph{fanfold}. A \emph{facet} of $\Sigma$ is a cone that is maximal for the inclusion. $\Sigma$ is \emph{pure dimensional} if all its facets have the same dimension. The \emph{$k$-skeleton of\/ $\Sigma$} is by definition the subfan of $\Sigma$ consisting of all the cones of dimension at most $k$, and we denote it by $\Sigma_{(k)}$. We say that $\Sigma$ is \emph{rational} if all of its cones are rational. $\Sigma$ is called \emph{simplicial} if each cone in $\Sigma$ is generated by as many rays as its dimension. It is called \emph{unimodular} if it is rational, simplicial, and the primitive vectors of the rays of any cone in $\Sigma$ are part of a basis for the entire lattice $N$.

The set of \emph{linear functions on $\Sigma$} is defined as the restriction to $\supp \Sigma$ of linear functions on $\ssN_{\R}$; such a linear function is defined by an element of $\ssM_{\R}$. In the case $\Sigma$ is rational, a linear function on $\Sigma$ is called \emph{integral} if it is defined by an element of $M$.

Let $f\colon \supp{\Sigma} \to \R$ be a continuous function. We say that $f$ is \emph{conewise linear on $\Sigma$} if on each face $\sigma$ of $\Sigma$, the restriction $f\rest\sigma$ of $f$ to $\sigma$ is linear. In such a case, we simply write $f\colon \Sigma \to \R$, and denote by $f_\sigma$ the linear form on $\ssN_{\sigma,\R}$ that coincides with $f\rest\sigma$ on $\sigma$. If the linear forms $f_\sigma$ are all integral, then we say $f$ is \emph{conewise integral linear}. (From an algebraic geometric point of view, these are meromorphic functions on $\Sigma$.)

A rational fan $\Sigma$ is called \emph{saturated at $\sigma$} in $\Sigma$ if the set of integral linear functions on $\ssSigma^\sigma$ coincides with the set of linear functions on $\ssSigma^\sigma$ that are conewise integral. This is equivalent to requiring the lattice generated by the points in $\supp{\ssSigma^\sigma}\cap \ssN^\sigma$ be saturated in the lattice $\ssN^\sigma$. We say $\Sigma$ is \emph{saturated} if it is saturated at each of its faces.

A conewise linear function $f \colon \Sigma \to \R$ is called \emph{convex}, \resp \emph{strictly convex}, if for each face $\sigma$ of $\Sigma$, there exists a linear function $\lambda$ on $\Sigma$ such that $f-\lambda$ vanishes on $\sigma$ and is non-negative, \resp strictly positive, on $\eta \setminus \sigma$ for any cone $\eta$ in $\Sigma$ that contains $\sigma$.

A fan $\Sigma$ is called \emph{quasi-projective} if it admits a strictly convex conewise linear function. When $\Sigma$ is rational, it is quasi-projective if and only if the toric variety $\P_\Sigma$ is quasi-projective.

The \emph{star fan} $\ssSigma^\sigma$ refers to the fan in $\ssN_{\R}^\sigma=\rquot {\ssN_{\R}}{\ssN_{\sigma,\R}}$ induced by the projection of the cones $\eta$ in $\Sigma$ that contain $\sigma$ as a face, for the projection map $\ssN_{\R} \to \ssN_{\R}^\sigma$.

\subsection*{Convention}

We endow the fan $\Sigma$ with the partial order $\subface$ given by the inclusion of cones in $\Sigma$: we write $\tau \subface \sigma$ if $\tau \subseteq \sigma$. We say $\sigma$ \emph{covers} $\tau$ and write $\tau \ssface \sigma$ if $\tau\subface\sigma$ and $\dims{\tau} = \dims\sigma-1$. The \emph{meet} operation $\wedge$ is defined as follows. For two cones $\sigma$ and $\delta$ of $\Sigma$, we set $\sigma \wedge \delta \coloneqq \sigma \cap \delta$. The set of cones in $\Sigma$ that contain both $\sigma$ and $\delta$ is either empty or has a minimal element $\eta \in \Sigma$. In the latter case, we say that $\eta$ is the \emph{join} of $\sigma$ and $\delta$ and denote it by $\sigma \vee \delta \coloneqq \eta$. We write $\sigma \sim \delta$ if $\sigma \vee \delta$ exists and $\sigma \wedge \delta = \conezero$ holds.

\subsection{Canonical compactification}

\begin{figure}[h]
  \centering
\begin{minipage}{0.45\textwidth}
  \centering
  \includegraphics{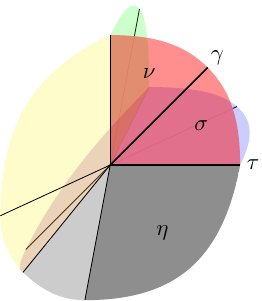}
\end{minipage}
\begin{minipage}{0.45\textwidth}
  \centering
  \includegraphics{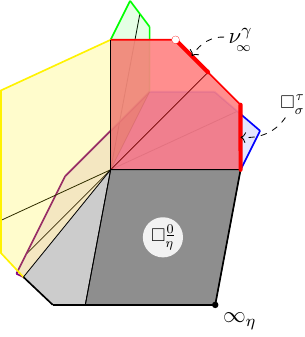}
\end{minipage}
\caption{A fan on the left and its canonical compactification on the right. The cone $\sigma$ is two-dimensional and has rays $\tau$ and $\gamma$. Two faces $\cube^\tau_\sigma$ and $\cube^\conezero_\eta$ are depicted in red and gray, respectively. The ray $\nu^\gamma_\infty$ depicted in bold red is based at the point $\ssinfty_\gamma$ and excludes the point $\ssinfty_{\oldnu}$, depicted in white.}
\label{fig:canonical_compactification}
\end{figure}

We briefly discuss canonical compactifications of fans and their combinatorics. More details can be found in \cites{AP-tht, OR11}.

Let $\Sigma$ be a fan in $\ssN_{\R}$. For any cone $\sigma$, denote by $\sssigma^\vee$ the \emph{dual cone} defined by
\[\sssigma^\vee \coloneqq \left\{m \in \ssM_{\R} \:\bigm|\: \langle m, a \rangle \geq 0 \:\textrm{ for all } a \in \sigma\right\}. \]

The \emph{canonical compactification} $\comp\sigma$ of $\sigma$ is given by $\hom_{\RMod}(\sssigma^\vee, \T_+)$, \ie, by the set of morphisms $\sssigma^\vee \to \T_+$ in the category of $\R_+$-modules. In this definition, we can naturally identify $\sigma$ with the corresponding subset of $\comp\sigma$. There is a natural topology on $\comp\sigma$ that makes it a compact topological space whose induced topology on $\sigma$ coincides with the Euclidean one. For an inclusion of cones $\tau \subseteq \sigma$, we get an inclusion map $\comp \tau \subseteq \comp \sigma$ that identifies $\comp \tau$ as the topological closure of $\tau$ in $\comp \sigma$.

The \emph{canonical compactification} $\comp\Sigma$ is defined as the union of $\comp\sigma$, $\sigma\in \Sigma$, where for an inclusion of cones $\tau \subseteq \sigma$ in $\Sigma$, we identify $\comp\tau$ with the corresponding subspace of $\comp\sigma$. An example of a canonical compactification is depicted in Figure~\ref{fig:canonical_compactification}. The topology of $\comp\Sigma$ is the induced quotient topology. Each extended cone $\comp\sigma$ naturally embeds as a subspace of $\comp\Sigma$.

There is a special point $\ssinfty_\sigma$ in $\comp\sigma$ defined by the map $\sssigma^\vee\to\T$ that vanishes on the orthogonal space $\sssigma^\perp \coloneqq \{m\in M_{\R}\:|\: \langle m, a \rangle = 0 \:\textrm{ for all } a \in \sigma\}$ and takes value $\ssinfty$ everywhere else. Note that for the cone $\conezero$, we have $\ssinfty_{\conezero} = 0$.

The compactification $\comp\Sigma$ naturally lives in the tropical toric variety $\ssTP_\Sigma$ defined as follows. For $\sigma\in \Sigma$, let $\chart\sigma \coloneqq \hom_{\RMod}(\sssigma^\vee,\T)$ and note that, since $\hom_{\RMod}(\sssigma^\vee, \R)\simeq \ssN_{\R}$, this is a partial compactification of $\ssN_{\R}$. For a pair of elements $\tau \subseteq \sigma$ in $\Sigma$, we get an inclusion $\~\tau\subseteq\~\sigma$. Gluing the spaces $\chart\sigma$ along these inclusions gives $\ssTP_\Sigma$.

We set $\ssN_{\infty,\R}^\sigma\coloneqq \ssN_{\R} + \ssinfty_\sigma \subseteq \~\sigma$. In this notation, we have $\ssN_{\infty,\R}^\conezero=\ssN_{\R}$. More generally, we have an isomorphism $\ssN_{\infty,\R}^\sigma \simeq \ssN_{\R}^\sigma$.

The tropical toric variety $\ssTP_\Sigma$ is naturally stratified as the disjoint union of tropical torus orbits $\ssN_{\infty,\R}^\sigma \simeq \ssN_{\R}^\sigma$, $\sigma \in \Sigma$. The natural inclusion of $\comp\sigma$ into $\chart\sigma$ gives an embedding $\comp\Sigma \subseteq \ssTP_\Sigma$ that identifies $\comp \Sigma$ as the closure of $\Sigma$ in $\ssTP_\Sigma$.

\subsection{Stratification of \texorpdfstring{$\comp\Sigma$}{the canonical compactification}} \label{sec:conical_stratification}

Consider a cone $\sigma \in \Sigma$ and a face $\tau$ of $\sigma$. Let $\sssigma^\tau_{\infty}$ be the subset of $\comp\sigma$ defined by
\[\sssigma^\tau_{\infty} \coloneqq \{\ssinfty_\tau + x \mid x\in \sigma\} = \comp\sigma \cap \ssN^\tau_{\infty,\R}. \]
Under the natural isomorphism $\ssN_{\infty,\R}^\tau \simeq \ssN_{\R}^\tau$, $\sssigma^\tau_{\infty}$ becomes isomorphic to the projection of the cone $\sigma$ into the linear space $\ssN_{\R}^\tau \simeq \ssN_{\infty,\R}^\tau$. We denote by $\opencone{\sigma}^\tau_{\infty}$ the relative interior of $\sssigma^\tau_{\infty}$. In order to simplify the notation, in what follows we will denote the closure $\sscompsigma^\tau_\infty$ of $\sssigma^\tau_{\infty}$ by $\cube^\tau_\sigma$. This is justified by observing that when $\Sigma$ is simplicial (assumption we make in the paper), $\sscompsigma^\tau_\infty$ is isomorphic to the hypercube $\T_+^k$, $k=\dims\sigma-\dims\tau$.

Note that the origin of $\sssigma^\tau_\infty$ is the point $\ssinfty_\tau$. By an abuse of the notation, we use $\ssinfty_\tau$ for the cone $\ss\tau_\infty^\tau =\{\ssinfty_\tau\}$.

The following proposition gives a precise description of how these different sets are positioned together in the canonical compactification.

\begin{prop}\label{prop:conical_stratification}
Let $\Sigma$ be a fan in $\ssN_{\R}$.
\begin{itemize}
\item The canonical compactification $\comp\Sigma$ is the disjoint union of (open) cones $\opencone{\sigma}^\tau_{\infty}$ for pairs of faces $\tau \subface \sigma$. The linear span of the cone $\sssigma^\tau_{\infty}$ is the real vector space $\ssN_{\infty, \sigma,\R}^\tau$, \ie, the projection of $\ssN_{\sigma, \R}$ into $\ssN_{\infty, \R}^\tau$.
\item For any pair of faces $\tau \subface \sigma$, the closure $\cube^\tau_\sigma$ of $\sssigma^\tau_{\infty}$ in $\comp\Sigma$ is the disjoint union of all the (open) cones $\opencone{\eta}^{\delta}_{\infty}$ with $\tau \subface \delta \subface \eta \subface \sigma$.
\end{itemize}
\end{prop}

\begin{proof}
The proof is a consequence of the tropical orbit-stratum correspondence theorem in the tropical toric variety $\TP_\Sigma$ and the observation we made previously that $\comp \Sigma$ is the closure of $\Sigma$ in $\TP_\Sigma$. We omit the details.
\end{proof}

The cones $\opencone{\sigma}^\tau_{\infty}$ form the \emph{open faces} of what we call the \emph{conical stratification} of $\comp\Sigma$. We refer to the topological closures $\cube^\tau_\sigma = \sscompsigma^\tau_{\infty}$ as the \emph{closed faces}, or simply, \emph{faces} of $\comp\Sigma$.

The closed faces $\cube^{\tau}_\sigma$, for $\tau\subface\sigma$ a pair of faces of $\Sigma$, endow $\comp\Sigma$ with an \emph{extended polyhedral structure}, see~\cites{JSS19, IKMZ} for the definition.

We extend to $\comp\Sigma$ the notations introduced for simplicial complexes. In particular, $\delta \in \comp\Sigma$ means that $\delta$ is a face of $\comp\Sigma$, $\supp{\comp\Sigma}$ denotes the support of $\comp\Sigma$, and $\compSigma_k$ is the set of faces of dimension $k$ in $\comp\Sigma$.

The \emph{sedentarity} of a face $\delta = \cube^\tau_\sigma$ denoted by $\sed(\delta)$ is by definition the face $\tau$ of $\Sigma$. By an abuse of the terminology and remembering only the dimension of $\sed(\delta)$, we sometimes say that $\delta$ has sedentarity $k$ with $k$ the dimension of $\tau$.

The tangent space to $\cube^\tau_\sigma$ is identified with $\ssN^{\tau}_{\infty,\sigma, \R}$, and contains full rank lattice $\ssN^{\tau}_{\infty, \sigma}$. For an inclusion of faces $\delta' = \cube^{\tau'}_{\sigma'} \subface \delta = \cube^{\tau}_{\sigma}$, we have $\tau \subface \tau' \subface \sigma' \subface \sigma$.

\subsection{Unit normal vectors and canonical multivectors and forms} \label{sec:canonical_forms}

Let $\Sigma$ be a rational fan of dimension $d$. Let $\sigma$ be a cone of $\Sigma$ and let $\tau$ be a face of codimension one in $\sigma$. Then, $\ssN_{\tau,\R}$ cuts $\ssN_{\sigma,\R}$ into two closed half-spaces only one of which contains $\sigma$. Denote this half-space by $\ssH_\sigma$. By a \emph{unit normal vector to $\tau$ in $\sigma$} we mean any vector $v$ of $\ssN_\sigma \cap \ssH_\sigma$ such that $\ssN_\tau + \Z v = \ssN_\sigma$. We usually denote such an element by $\nvect_{\sigma/\tau}$ and note that it induces a well-defined generator of $\ssN^\tau_\sigma = \rquot{\ssN_\sigma}{\ssN_\tau}$ that we denote by $\sse^\tau_\sigma$. We naturally extend the definition to similar pair of faces in $\comp\Sigma$ having the same sedentarity. In the case $\sigma$ is a ray (and $\tau$ is a point of the same sedentarity as $\sigma$), we also use the notation $\sse_\sigma$ instead of $\nvect_{\sigma/\tau}$.

On each face $\sigma$ of $\Sigma$, we fix a generator of $\bigwedge^{\dims{\sigma}}\ssN_\sigma$ denoted by $\ssnu_\sigma$ and call it the \emph{canonical multivector of $\sigma$} (up to a sign, this is unique). The element $\ssvarpi_\sigma\in \bigwedge^{\dims{\sigma}}\ssN_\sigma^\dual$ that takes value one on $\ssnu_\sigma$ is called the \emph{canonical form of $\sigma$}. We assume that $\ssnu_{\conezero} = 1$ and that $\ssnu_\rho = \sse_\rho$ for any ray $\rho$.

We extend the above definitions to $\comp\Sigma$ by choosing an element $\ssnu^\tau_\sigma \in \bigwedge^{\dims{\sigma}-\dims{\tau}}\ssN^\tau_{\infty,\sigma}\simeq \bigwedge^{\dims{\sigma}-\dims{\tau}}\ssN^\tau_\sigma$ for each polyhedron $\cube^\tau_\sigma$, and setting $\varpi^\tau_\sigma \coloneqq \ssnu^{\tau\,\dual}_\sigma$. For simplicial $\Sigma$, there is a natural choice for extending the cellular orientation of $\Sigma$ to $\comp\Sigma$. If $\tau \subface \sigma$ is a pair of cones, there exists a unique minimal $\tau'$ such that $\sigma = \tau \vee \tau'$. Then, we set $\ssnu^\tau_\sigma$ to be the image of $\ssnu_{\tau'}$ via the projection $\ssN \to \ssN_{\infty, \R}^\tau\simeq \ssN^\tau$.

\subsection{Orientation and sign function} \label{subsec:orientation}

The choice of canonical multivectors from the previous section associate an orientation (in the topological sense) to the faces of $\comp\Sigma$. In particular, to any pair of closed faces $\gamma \ssface \delta$ of $\comp\Sigma$, we associate a \emph{sign} denoted $\sign(\gamma,\delta)$ as follows.

If both faces have the same sedentarity, then $\sign(\gamma,\delta)$ is the sign of $\ssvarpi_\delta\bigl(\nvect_{\delta/\gamma} \wedge \ssnu_\gamma\bigr)$. Otherwise, there exists a pair of cones $\tau' \ssface \tau$ and a cone $\sigma$ in $\Sigma$ such that $\gamma = \cube^\tau_\sigma$ and $\delta = \cube^{\tau'}_\sigma$. Consider the map $\sspi_{\gamma \ssface \delta} \colon \ssN_{\infty, \sigma}^{\tau'} \to \ssN_{\infty, \sigma}^\tau$ between the two lattices in the tangent spaces of $\gamma$ and $\delta$ (see Section~\ref{sec:conical_stratification}). Note that $\sspi_{\gamma \ssface \delta}$ is surjective and induces a surjective linear map from $\bigwedge^{\bul}\ssN^{\tau'}_{\infty,\sigma}$ to $\bigwedge^{\bul}\ssN^{\tau}_{\infty,\sigma}$. We choose $\ssnu'$ such that $\sspi_{\gamma\ssface\delta} (\ssnu') = \ssnu_{\gamma}$. Let $\sse^{\tau'}_\tau$ be the primitive vector of the ray $\sstau^{\tau'}_{\infty}$ in $\sssigma^{\tau'}_{\infty}$. Although, there is a choice for $\ssnu'$, the element $\sse^{\tau'}_\tau \wedge \ssnu'$ is well-defined and does not depend on this choice. We define $\sign(\gamma,\delta)$ as the sign of $-\ssvarpi_\delta(\sse^{\tau'}_\tau \wedge \ssnu')$.

\subsection{Chow rings of fans and localization lemma}
\label{sec:localization_lemma}

The \emph{Chow ring} of a rational simplicial fan $\Sigma$ denoted by $A^\bul(\Sigma)$ is defined by generators and relations as follows. Consider the polynomial ring $\Z[\x_\zeta]_{\zeta\in \ssSigma_1}$ with indeterminate variables $\x_\zeta$ associated to rays $\zeta$ in $\ssSigma_1$. The ring $A^\bul(\Sigma)$ is the quotient ring
\[ A^\bul(\Sigma) \coloneqq \rquot{\Z[\x_\zeta]_{\zeta\in \Sigma_1}}{\bigl(I + J\bigr)} \]
where $I$ is the ideal generated by the products $\x_{\rho_1}\!\cdots \x_{\rho_k}$, for $k\in \N$, such that $\ssrho_1, \dots, \ssrho_k$ are non-comparable rays in $\Sigma$, that is, they do not form a cone in $\Sigma$, and $J$ is the ideal generated by the elements of the form
\[ \sum_{\zeta\in \Sigma_1} m(\sse_\zeta)\x_\zeta, \qquad  m \in M \coloneqq \ssN^\dual,\]
with $\sse_\zeta$ the primitive vector of the ray $\zeta$.

The ideal $I+J$ is homogeneous and the Chow ring inherits a graded ring structure. Denoting by $d$ the dimension of $\Sigma$, the graded pieces for degree larger than $d$ vanish (by Theorem~\ref{thm:localization-lemma} below), and we can write
\[ A^\bul(\Sigma) = \bigoplus_{k= 0}^d A^k(\Sigma). \]
We define the Chow ring with rational and real coefficients in a similar way and denote them by $A^\bul(\Sigma, \Q)$ and $A^\bul(\Sigma, \R)$, respectively. A useful result dealing with the Chow rings is the localization lemma stated below that provides a reformulation of each graded piece $A^k(\Sigma)$ of the Chow ring.

For each cone $\sigma$ of dimension $k$ in $\Sigma$, choose an indeterminate variable $\x_\sigma$ and denote by $Z^k(\Sigma)$ the free abelian group generated by $\x_\sigma$ for $\sigma\in \Sigma_k$. Consider the injective map $Z^k(\Sigma) \hookrightarrow \Z[\x_\rho \mid \rho\in \Sigma_1]$ that sends each $\x_\sigma$ to the product $\prod_{\rho}\x_\rho$ over rays $\rho$ of $\sigma$. This induces a map of abelian groups $Z^k(\Sigma) \to A^k(\Sigma)$.

We have the following result proved in~\cite[Theorem 1]{FMSS}, see~\cite[Theorem 3.2]{AP-hodge-fan} for a combinatorial proof.

\begin{thm}[Localization lemma] \label{thm:localization-lemma}
  Notations as above, if\/ $\Sigma$ is unimodular, the map from $Z^k(\Sigma)$ to $A^k(\Sigma)$ is surjective and its kernel is generated by elements of the form $\sum_{\sigma\in\Sigma_k \\
\tau\ssface \sigma} m(\sse_{\sigma}^{\tau})\x_\sigma$ for $\tau \in \Sigma_{k-1}$ and $m \in M^\tau=(N^\tau)^\dual$.

The same statement holds true with rational coefficients for a simplicial rational fan.
\end{thm}

\subsection{Tropical homology and cohomology groups} \label{subsec:homology}

Let $\Sigma$ be a rational fan. The extended polyhedral structure on $\comp\Sigma$ leads to the definition of tropical homology groups and cohomology groups introduced in~\cite{IKMZ} and further studied in~\cites{JSS19, MZ14, JRS18, GS-sheaf, AB14, ARS, Aks21}.

We recall the definition of the multi-tangent and multi-cotangent (integral) spaces $\SF_p$ and~$\SF^p$. We give the definitions for $\comp\Sigma$ for the face structure given by the closed strata $\cube^\tau_\sigma$ in the conical stratification of $\comp\Sigma$, $\tau \subface \sigma$ in $\Sigma$. This gives a combinatorial complex which calculates the tropical homology and cohomology groups of $\comp \Sigma$. The definitions are naturally adapted to $\Sigma$ with face structure given by its cones.

For any face $\delta = \cube^\tau_\sigma$, we set $N_\delta = \ssN_{\infty, \sigma}^\tau \simeq \rquot{\ssN_{\sigma}}{\ssN_{\tau}}$. For any non-negative integer $p$, the \emph{$p$-th multi-tangent} and the \emph{$p$-th multi-cotangent space of\/ $\comp\Sigma$ at $\delta$}, denoted by $\SF_p(\delta)$ and $\SF^p(\delta)$ respectively, are given by
\[\SF_p(\delta)\coloneqq \hspace{-.3cm} \sum_{\eta \supface \delta \\ \sed(\eta) = \tau } \hspace{-.3cm} \bigwedge^p N_\eta \ \subseteq \bigwedge^p N^\tau, \quad \textrm{and} \quad \SF^p(\delta) \coloneqq \SF_p(\delta)^\dual.\]

For an inclusion of faces $\gamma \subface \delta$ in $\comp\Sigma$, we get maps $\ssi_{\delta \supface \gamma}\colon \SF_p(\delta) \to \SF_p(\gamma)$ and $\ssi_{\gamma \subface \delta}^*\colon \SF^p(\gamma) \to \SF^p(\delta)$ defined as follows. If $\gamma$ and $\delta$ have the same sedentarity, the map $\ssi_{\delta \supface \gamma}$ is just an inclusion. If $\gamma=\cube_\sigma^{\tau'}$ and $\delta=\cube_\sigma^\tau$ with $\tau \subface \tau' \subface \sigma$, then the map $\ssi_{\delta\supface \gamma}$ is induced by the projection $N^\tau_{\infty} \to N^{\tau'}_{\infty}$. In the general case, $\ssi_{\delta\supface \gamma}$ is given by the composition of the projection and the inclusion; the map $\ssi_{\gamma \subface \delta}^*$ is the dual of $\ssi_{\delta\supface \gamma}$.

Let $X = \Sigma$ or $\comp\Sigma$, for a fan $\Sigma$. For a pair of non-negative integers $p,q$, define
\[ C_{p,q}(X) \coloneqq \bigoplus_{\delta \in X \\ \dims\delta=q \\ \delta\text{ compact}} \SF_p(\delta) \]
and consider the corresponding complexes
\[C_{p,\bullet}(X)\colon \quad \dots\longrightarrow C_{p, q+1}(X) \xrightarrow{\ \partial_{q+1}\ } C_{p,q}(X) \xrightarrow{\ \partial_q\ } C_{p,q-1} (X)\longrightarrow\cdots\]
where the differential is given by the sum of maps $\sign(\gamma,\delta)\cdot\ssi_{\delta\ssupface \gamma}$ with the signs corresponding to a chosen cellular orientation on $X$ as explained in Section~\ref{subsec:orientation}.

\smallskip
The \emph{tropical homology} with integral coefficients of $X$ is defined by
\[ H_{p,q}(X) \coloneqq H_q(C_{p,\bullet}(X)). \]

Similarly, we have a cochain complex
\[C^{p,\bullet}(X)\colon \quad \dots\longrightarrow C^{p, q-1}(X) \xrightarrow{\ \d^{q-1}\ } C^{p,q}(X) \xrightarrow{\ \d^{q}\ } C^{p,q+1}(X) \longrightarrow \cdots\]
where
\[C^{p,q}(X) \coloneqq C_{p,q}(X)^\dual \simeq \! \bigoplus_{\delta \in X \\ \dims\delta = q \\ \delta \text{ compact}} \!\! \SF^p(\delta)\]
and the \emph{tropical cohomology} with integral coefficients of $X$ is defined by
\[H^{p,q}(X) \coloneqq H^q(C^{p,\bullet}(X)).\]

We can also define the \emph{compact-dual versions} of tropical homology and cohomology by allowing non-compact faces. These are called Borel-Moore homology and cohomology with compact support, and are defined as follows.
\[ C^\BM_{p,q}(X) \coloneqq \bigoplus_{\delta \in X \\ \dims\delta=q} \SF_p(\delta) \quad \text{and} \quad C_c^{p,q}(X) \coloneqq \bigoplus_{\delta \in X \\ \dims\delta = q } \SF^p(\delta). \]
We get the corresponding (co)chain complexes $C^\BM_{p,\bul}(X)$ and $C_c^{p,\bul}(X)$, and the \emph{Borel-Moore tropical homology} and the \emph{tropical cohomology with compact support} with integral coefficients are respectively
\[ H^\BM_{p,q}(X) \coloneqq H_q(C^\BM_{p,\bul}(X)) \quad \text{and} \quad H_c^{p,q}(X) \coloneqq H^q(C_c^{p,\bul}(X)). \]
If $X$ is compact, then both notions of homology and both notions of cohomology coincide.

Similarly, we define homology and cohomology groups with rational coefficients.

Homology and cohomology in this paper refer to the tropical ones, so we usually omit the mention of the word tropical. In our definition of tropical homology and cohomology, we adapted a cellular version. As in the classical setting, there exist other ways of computing the same groups: for instance using either of singular, cubical, or sheaf cohomologies. We note in particular that the homology and cohomology only depends on the support.

We need the following result on the cohomology with compact support of $\T^k$, see~\cite{JSS19} for the proof.

\begin{prop}\label{prop:cohomology_compact_support_Tk} The cohomology $H_c^{p,q}(\T^k)$ with compact support of\/ $\T^k$ is trivial unless $p=q=k$, and we have $H_c^{p,q}(\T^k) \simeq \Z$.
\end{prop}

The collection of coefficient sheaves and cosheaves come with contraction maps defined as follows.

Let $\sigma \in \Sigma$ and consider an element $\ssnu \in \SF_{k}(\sigma)$. Given an element $\alpha \in \SF^{p}(\sigma)$, we denote by $\ss\contract_{\ssnu}(\alpha)$ the element of $\SF^{p-k}(\sigma)$ that on each element $\ssnu' \in \SF_{p-k}(\sigma)$ takes value $\alpha(\ssnu \wedge \ssnu')$. The map $\ss\contract_{\ssnu} \colon \SF^{p}(\sigma) \to \SF^{p-k}(\sigma)$ is linear. We extend the definition naturally to all faces $\delta$ of $\comp\Sigma$ and obtain maps $\ss\contract_{\ssnu} \colon \SF^{p}(\delta) \to \SF^{p-k}(\delta)$ for any element $\ssnu \in \SF_k(\delta)$. We call $\ss\contract_{\ssnu}(\alpha)$ the contraction of $\alpha$ by $\ssnu$, and refer to $\ss\contract_{\ssnu}$ as the contraction map defined by $\ssnu$. Dually, for $\alpha \in \SF^{k}(\sigma)$, we get contraction maps $\contract_\alpha \colon \SF_p(\sigma) \to \SF_{p-k}(\sigma)$.

These contraction maps lead to the definition of cap product between homology and cohomology groups, see next section for cap product with the fundamental class.

\subsection{The case of tropical fans} \label{sec:tropical-fans}

A \emph{tropical orientation} of a rational fan $\Sigma$ of pure dimension $d$ is an integer valued map
\[ \omega\colon \Sigma_d \to \Z\setminus\{0\} \]
that verifies the \emph{balancing condition}. This means that for any cone $\tau$ in $\Sigma$ of codimension one, the sum
\[ \sum_{\sigma \ssupface \tau} \omega(\sigma)\e_{\sigma}^{\tau} \]
vanishes in the quotient lattice $\rquot{N}{N_\tau}$.

A \emph{tropical fan} in $N$ is a pair $(\Sigma, \ssomega_\Sigma)$ consisting of a pure dimensional rational fan $\Sigma$ and a tropical orientation $\ssomega_\Sigma$ as above. We will call \emph{tropical fanfold} the support of any tropical fan $(\Sigma, \ssomega_\Sigma)$ endowed with the tropical orientation induced by $\ssomega_\Sigma$ on the set of regular points (those points having an open neighborhood isomorphic to an open subset of a real vector space).

In the case of a tropical fan $\Sigma$, the balancing condition implies the existence of a \emph{fundamental class} $[\Sigma]\in H_{d,d}^\BM(\Sigma)$ defined using canonical multivectors from Section~\ref{sec:canonical_forms}. The class $[\Sigma]$ is represented by \emph{the canonical element} $\ssnu_{\Sigma} \in C_{d,d}^{\BM}$ given by
\[ \ssnu_{\Sigma} \coloneqq \bigl(\ssomega_\Sigma(\eta)\ssnu_\eta\bigr)_{\eta \in \Sigma_d} \in \bigoplus_{\eta \in \Sigma_d} \bigwedge^d \ssN_\eta. \]

Cap product with the fundamental class induces a map $\frown [\Sigma] \colon H^{p,q}(\Sigma)\to H^{\BM}_{d-p,d-q}(\Sigma)$ for each $p,q\in \{0,\dots,d\}$. We say that $\Sigma$ satisfies \emph{tropical Poincaré duality with integral coefficients} if these maps are all isomorphisms for all $p$ and $q$. A tropical fanfold $\supp{\Sigma}$ is called a \emph{tropical homology manifold} if any open subset $U$ of $\supp{\Sigma}$ satisfies tropical Poincaré duality, \ie, cap product induces an isomorphism between the tropical cohomology and the tropical Borel-Moore homology of $U$ (see \cite{JSS19,JRS18} for details). Equivalently, a tropical fanfold $\supp{\Sigma}$ is a tropical homology manifold if any star fan $\ssSigma^\sigma$ of $\Sigma$ verifies tropical Poincaré duality.

Note that $H^{p,q}(\Sigma)$ is trivial for $q \neq 0$, and for $q=0$, we have $H^{p,q}(\Sigma) =\SF^p(\conezero)$. Cap product with fundamental class is described using contraction maps defined in the previous section. In the case $q=0$, it is the map
\[ \SF^k(\conezero) \to H^\BM_{d-k, d}(\Sigma) \]
defined as follows. Consider the fundamental class $[\Sigma] \in H^\BM_{d,d}(\Sigma)$. The coefficient of a facet $\sigma\in \Sigma_d$ in $[\Sigma]$ is given by $\omega(\sigma)\nu_\sigma$. Given an element $\alpha \in \SF^k(\conezero)$, the corresponding element in $H^\BM_{d-k, d}(\Sigma)$ has coefficient at the facet $\sigma \in \Sigma_d$ given by $\contract_{\alpha}(\omega(\sigma)\nu_{\sigma})$. Poincaré duality for $\Sigma$ is the statement that these maps are all isomorphisms. Note that, dually, we get a map
\[ H_c^{d-k, d}(\Sigma) \to \SF_k(\conezero). \]

\subsection{Minkowski weights} \label{sec:MW}

Let $\Sigma$ be a rational fan in $N$. A \emph{Minkowski weight} of dimension $p$ on $\Sigma$ with coefficient in $\Z$ is a map $w \colon \Sigma_{p} \to \Z$ that verifies the balancing condition, namely, that
\[ \forall\, \tau\in \Sigma_{p-1},\qquad \sum_{\sigma\ssupface\tau} w(\sigma) \sse^\tau_\sigma = 0 \in \ssN^{\tau}. \]

We denote by $\MW_{p}(\Sigma)$ the set of all Minkowski weights of dimension $p$ on $\Sigma$ with integral coefficients. Addition of weights cell by cell turns $\MW_{p}(\Sigma)$ into a group.

Replacing $\Z$ with $\Q$, we get the set of Minkowski weights with rational coefficients. It coincides with the vector space generated by $\MW_{p}(\Sigma)$, and is denoted by $\MW_{p}(\Sigma, \Q)$.

The \emph{support} of a Minkowski weight $w$ of dimension $p$ is the set of $\sigma \in \Sigma_p$ with $w(\sigma) \neq 0$.

Given a Minkowski weight $w$ of dimension $p$, we consider its support and its closure in $\comp\Sigma$ which is a tropical cycle of dimension $p$. Taking the corresponding homology class in $ H_{p,p}(\comp\Sigma)$ leads to the cycle class map $\cl\colon \MW_p(X) \to H_{p,p}(\comp\Sigma)$, see for example~\cites{AP20-hc, MZ14, Sha-thesis, JRS18, GS19}. As a consequence of the localization lemma, we have the following duality between Minkowski weights and Chow groups~\cite{AHK,AP-hodge-fan}.
\begin{thm}\label{thm:duality_chow_mw} Let $\Sigma$ be a unimodular fan. There is a pairing
  \begin{align*}
    A^k(\Sigma) \otimes \MW_k(\Sigma) &\to \Z\\
    \x_\sigma\, , w \ \quad\qquad &\mapsto w(\sigma) \quad \quad \forall \sigma \in \Sigma_k \textrm{ and } w \in\MW_k(\Sigma).
  \end{align*}
 This induces an isomorphism
 \[ \MW_k(\Sigma) \simeq A^k(\Sigma)^*. \]
With rational coefficients, the same statement holds true for simplicial rational fans.
\end{thm}

Note that in general $A^k(\Sigma)$ can have torsion, and the pairing in the theorem is not necessarily perfect.

\section{The fine double complex and proof of Theorem~\ref{thm:cohomology_compactification}} \label{sec:fine_double_complex}

We start with the proof of Theorem~\ref{thm:cohomology_compactification}. Following our convention, we give the proof of the statement assuming $\Sigma$ is a unimodular fan in $N_{\R}$, working with integral coefficients. The proof for simplicial rational fans with rational coefficients will be similar.

We fix a non-negative integer $p$.

\subsection{The fine double complex}

The tropical cochain complex $C^{p,\bul}(\comp\Sigma)$ given in Section~\ref{subsec:homology} is the total complex of a double complex. Namely, for each face $\cube^\tau_\sigma = \compcone{\sigma}^\tau_\infty$ of $\comp\Sigma$, we can remember both numbers $\dims\sigma$ and $\dims\tau$, instead of the grading by dimension of the face, $\dims\sigma-\dims\tau$, that resulted in the cochain complex $C^{p,\bul}(\comp\Sigma)$. Unfolding in this way $C^{p,\bul}(\comp\Sigma)$ leads to the double complex
\[ \Ep{p}^{a,b} \coloneqq \bigoplus_{\tau\subface\sigma \textrm{ faces of }\Sigma\\ \dims\sigma = a \\ \dims\tau = -b} \SF^p(\cube^\tau_\sigma). \]
Differentials $\Ep{p}^{a,b} \to \Ep{p}^{a+1,b}$ and $\Ep{p}^{a,b} \to \Ep{p}^{a,b+1}$ are given by maps between the coefficient groups $\SF^p$ of faces of $\comp\Sigma$ induced by the inclusions of faces $\cube^\tau_\sigma \hookrightarrow \cube^\tau_\eta$ with $\sigma\ssface\eta$, for the first differential, and inclusions of faces $\cube^\tau_\sigma \hookrightarrow \cube^\gamma_\sigma$ for $\gamma \ssface \tau$, for the second, respectively. We have
\[ C^{p,\bul}(\comp\Sigma) = \Tot^\bul(\Ep{p}^{\bul, \bul}). \]

In order to prove Theorem~\ref{thm:cohomology_compactification}, we calculate the cohomology of $C^{p,\bul}(\comp\Sigma)$ by using the spectral sequence associated to the filtration given by the columns of $\Ep{p}^{\bul, \bul}$. (The spectral sequence associated to the filtration by the rows will be used in the proof of Theorem~\ref{thm:smoothness_alternate} in Section~\ref{sec:proof-smoothness-criterion}.)

We denote by $\EpI{p}^{\bul,\bul}_0$ the $0$-th page of this spectral sequence which has abutment
\begin{equation}\label{sec:spectral_sequence_columns}
\EpI{p}^{\bul,\bul} \Longrightarrow H^{p,\bul}(\comp\Sigma).
\end{equation}

\subsection{Computation of the first page}

The proof of Theorem~\ref{thm:cohomology_compactification} is based on the following result. For non-negative integers $p,q$, define
\[ \cubC^{p,q}(\comp\Sigma) \coloneqq \bigoplus_{\sigma \in \Sigma \\ \dims\sigma = q} \SF^{p-q}(\ssinfty_\sigma),\]
with $\SF^{k}(\ssinfty_\sigma)$ the coefficient group used in Section~\ref{sec:prel} in the definition of tropical cohomology. Here, and in what follows, by slight abuse of notation, we use $\ssinfty_\sigma$ when referring to the face $\sssigma^{\sigma}_\infty=\{\ssinfty_\sigma\}$.

Denote by $\EpI{p}_1^{\bul,\bul}$ the first page of the spectral sequence.

\begin{prop} \label{prop:degeneration_spectral_sequence}
Notations as above, we have
\[ \EpI{p}_1^{a,b} =\begin{cases}
  \cubC^{p,a}(\comp\Sigma) & \text{for $b = 0$,} \\
  0 & \text{otherwise.}
\end{cases} \]
The differential on the first page for $b=0$ coincides with the differential $\cubd^\bul$ of $\cubC^{p,\bul}(\comp\Sigma)$.
\end{prop}

In order to prove this proposition, we decompose $\EpI{p}_1^{\bullet,\bullet}$ to a direct sum of subcomplexes, as follows. Fix an integer $a$. The $a$-th column in zeroth page $\EpI{p}_0^{a,\bul}$ of the spectral sequence is
\[ \EpI{p}_0^{a, \bul}\colon \quad \dots \longrightarrow \bigoplus_{\tau\subface\sigma \textrm{ faces of }\Sigma \\ \dims\sigma = a \\ \dims\tau = -b+1} \SF^p(\cube^\tau_\sigma) \longrightarrow \bigoplus_{\tau\subface\sigma \textrm{ faces of }\Sigma \\ \dims\sigma = a \\ \dims\tau = -b} \SF^p(\cube^\tau_\sigma) \longrightarrow \bigoplus_{\tau\subface\sigma \textrm{ faces of }\Sigma \\ \dims\sigma = a \\ \dims\tau = -b-1} \SF^p(\cube^\tau_\sigma) \longrightarrow \cdots \]
We decompose this complex as a direct sum of complexes $\EpI{p}_0^{\sigma,\bul}$ associated to each face $\sigma \in \Sigma_a$
\[ \EpI{p}_0^{\sigma,\bul} \,\, :\quad \dots \longrightarrow \bigoplus_{\tau\subface\sigma \\ \dims\tau = -b+1} \SF^p(\cube^\tau_\sigma) \longrightarrow \bigoplus_{\tau\subface\sigma \\ \dims\tau = -b} \SF^p(\cube^\tau_\sigma) \longrightarrow \bigoplus_{\tau\subface\sigma \\ \dims\tau = -b-1} \SF^p(\cube^\tau_\sigma) \longrightarrow \cdots \]
We thus have
\begin{equation}\label{eq:decomposition}
\EpI{p}_0^{a,\bul} = \bigoplus_{\sigma \in \Sigma_a} \EpI{p}_0^{\sigma,\bul}.
\end{equation}
Let $\EpI{p}_1^{\sigma,\bul}$ be the cohomology of $\EpI{p}_0^{\sigma,\bul}$, so that we have
\begin{equation}\label{eq:decomposition1}
\EpI{p}_1^{a,\bul} = \bigoplus_{\sigma \in \Sigma_a} \EpI{p}_1^{\sigma,\bul}.
\end{equation}

These cohomology groups are given by the following lemma.
\begin{lemma}[Hypercube vanishing lemma] \label{lem:hypercube_vanishing}
For any face $\sigma$ of\/ $\Sigma$, the cohomology $\EpI{p}_1^{\sigma,\bul}$ of the complex $\EpI{p}_0^{\sigma,\bul}$ is given by
\[ \EpI{p}_1^{\sigma,b} \simeq \begin{cases}
  \SF^{p-\dims\sigma}(\ssinfty_\sigma) & \text{if $b=0$,} \\
  0 & \text{otherwise.}
\end{cases} \]
This isomorphism is moreover induced by the natural map
\[ \begin{array}{ccccc}
  \EpI{p}_0^{\sigma,0} &\simeq& \SF^p(\comp\sigma) &\longrightarrow& \SF^{p-\dims\sigma}(\ssinfty_\sigma), \\
  && \alpha & \longmapsto & \sspi^\sigma_*(\contract_{\oldnu_\sigma}(\alpha))
\end{array} \]
with $\pi^\sigma$ referring to the projection $N\to \ssN^\sigma$, extended to the exterior algebra and restricted to the corresponding subspaces $\SF^p$. The contraction map $\ss\contract_{\ssnu}$ is defined in Section~\ref{subsec:homology}.
\end{lemma}

We postpone the proof of this lemma to Section~\ref{sec:proof-vanishing-lemma}.

\begin{proof}[Proof of Proposition~\ref{prop:degeneration_spectral_sequence}]
Using the hypercube vanishing lemma, and in view of the decomposition given in~\eqref{eq:decomposition1}, the first page of the spectral sequence is concentrated in the 0-th row and is given by
\[ \EpI{p}_1^{a,0} = \bigoplus_{\sigma\in \Sigma\\ \dims\sigma = a} \SF^{p-a}(\ssinfty_\sigma), \]
which is precisely $\cubC^{p,a}(\comp\Sigma)$.

It remains to check that the differentials coincide. Let $\tau\ssface\sigma$ be two faces and denote by $\d_{\tau\ssface\sigma}\colon \SF^{p-\dims\tau}(\ssinfty_\tau) \to \SF^{p-\dims\sigma}(\ssinfty_\sigma)$ the corresponding part of the differential in $\EpI{p}^{\bul,0}$. By Lemma~\ref{lem:hypercube_vanishing}, we have the following commutative diagram.
\[ \begin{tikzcd}
\SF^p(\comp\tau) \rar{\sign(\tau,\sigma)\i^*}\dar{\pi^\tau_*\circ\,\contract_{\ssnu_\tau}}& \SF^p(\comp\sigma) \dar{\pi^\sigma_*\circ\,\contract_{\ssnu_\sigma}} \\
\SF^{p-\dims\tau}(\ssinfty_\tau) \rar{\d_{\tau\ssface\sigma}}& \SF^{p-\dims\sigma}(\ssinfty_\sigma)
\end{tikzcd} \]
We infer that the bottom map should be given by $\sspi_*\circ\ss\contract_{\e^\tau_\sigma}$ with $\pi\colon \ssN^\tau \to \ssN^\sigma$, \ie, by the differential of $\cubC^{p,\bul}(\comp\Sigma)$. This concludes the proof.
\end{proof}

\subsection{Proof of Theorem~\ref{thm:cohomology_compactification}}

This is a direct consequence of Proposition~\ref{prop:degeneration_spectral_sequence}, which shows that the spectral sequence degenerates at page two, and, moreover, we have
\[ H^{p,q}(\comp \Sigma) \simeq \Ep{p}_2^{q,0} = H^q\bigl(\cubC^{p,\bul}(\comp\Sigma)\bigr). \pushQED{\qed}\qedhere \]

\subsection{Toric weight filtration on $\SF^{\bul}$}

It remains to prove Lemma~\ref{lem:hypercube_vanishing}. For this, we introduce a natural filtration on $\SF^{\bul}$ called \emph{toric weight filtration}. We only use here basic properties of this filtration. We note however that this plays an important role in the development of Hodge theory for tropical varieties in our work~\cite{AP-tht}, and we refer to \emph{loc.cit.} for a more through study.

For each pair of faces $\tau \subface \sigma$ and each integer $k \leq p$, we have a map
\[ \begin{array}{ccl}
  \SF^p(\cube^\tau_\sigma) & \longrightarrow & \hom\bigl(\bigwedge^{k+1} \ssN^{\tau}_\sigma, \SF^{p-k-1}(\cube^\tau_\sigma)\bigr) \simeq \bigwedge^{k+1} \ssN^{\tau\,\dual}_\sigma \otimes \SF^{p-k-1}(\cube^\tau_\sigma), \\
  \alpha & \longmapsto & \hspace{10.5ex}(\nu \mapsto \ss\contract_{\oldnu}(\alpha)).
\end{array} \]
We denote by $\ssW_k\SF^p(\cube^\tau_\sigma)$ the kernel of this map. We get a filtration
\[ 0=\ssW_{-1}(\SF^p(\cube^\tau_\sigma)) \subseteq \ssW_{0}(\SF^p(\cube^\tau_\sigma)) \subseteq \cdots \subseteq \ssW_p(\SF^p(\cube^\tau_\sigma)) = \SF^p(\cube^\tau_\sigma) \]
of $\SF^p(\cube^\tau_\sigma)$. Let $\alpha \in \ssW_k\SF^p(\cube^\tau_\sigma)$ and $\nu \in \bigwedge^k \ssN^\tau_\sigma$, and let $\pi\colon \ssN^\tau \to \ssN^\sigma$ be the projection. Then, $\sspi_*(\ss\contract_\oldnu(\alpha)) \in \SF^{p-k}(\ssinfty_\sigma)$ is a well-defined element. This way, we obtain a map
\[ \ssW_k\SF^p(\cube^\tau_\sigma) \to \bigwedge^k \ssN_\sigma^{\tau\,*} \otimes \SF^{p-k}(\infty_\sigma). \]
It is straightforward to check that this map is surjective, and that its kernel is precisely $\ssW_{k-1}\SF^p(\cube^\tau_\sigma)$. We thus get the following description of the graded pieces of the filtration
\[ \gr_k(\ssW_\bul\SF^p(\cube^\tau_\sigma)) \coloneqq \rquot{\ssW_k\SF^p(\cube_\sigma^\tau)}{\ssW_{k-1}\SF^p(\cube_\sigma^\tau)} \simeq \bigwedge^k \ssN^{\tau\,\dual}_\sigma \otimes \SF^{p-k}(\ssinfty_\sigma). \]

\subsection{Proof of the hypercube vanishing lemma}\label{sec:proof-vanishing-lemma}

We calculate the cohomology of $\EpI{p}^{\sigma,\bul}_0$ via the spectral sequence $\Fpq{\sigma}{p}^{\bul,\bul}$ associated to the filtration $\ssW_{\bul}$. By the discussion of the previous section, the zeroth page of this spectral sequence is given by
\[ \Fpq{\sigma}{p}_0^{a,b} = \bigoplus_{\tau\subface \sigma \\ \dims\tau=-a-b} \bigwedge^a \ssN_\sigma^{\tau\,\dual} \otimes \SF^{p-a}(\ssinfty_\sigma). \]
The $a$-th row of the first page of this spectral sequence is given by the cohomology of the cochain complex $\Fpq{\sigma}{p}_1^{a,\bul}$
\[\Fpq{\sigma}{p}_1^{a,b} = H^b(\Fpq{\sigma}{p}_1^{a,\bul}).\]
To compute this cohomology, we notice that
\[ \Fpq{\sigma}{p}^{a,\bul}_0 \simeq \Bigl( \dots
  \longrightarrow \bigoplus_{\tau\subface\sigma \\ \dims\tau = -a-b+1} \bigwedge^{a}\ssN^{\tau\,\dual}_\sigma
  \longrightarrow \bigoplus_{\tau\subface\sigma  \\ \dims\tau = -a-b} \bigwedge^{a}\ssN^{\tau\,\dual}_\sigma
\longrightarrow \cdots \Bigr)\ \otimes \SF^{p-a}(\ssinfty_\sigma). \]
Since $\sigma$ is unimodular, the cochain complex appearing in the above summand is isomorphic to the cochain complex $C_c^{a,\bul+\dims\sigma+a}(\T^{\dims\sigma})$ for the cohomology with compact support and integral coefficients of $\T^{\dims\sigma}$ for the cell decomposition of $\T^{\dims\sigma}$ given by $\{\T, \ssinfty\}^{\dims\sigma}$:
\[ \Fpq{\sigma}{p}^{a,\bul}_0 \simeq C_c^{a,\bul+\dims\sigma+a}(\T^{\dims\sigma}) \otimes \SF^{p-a}(\ssinfty_\sigma). \]
(Without the unimodularity assumption on $\sigma$, this still holds with rational coefficients.)

We deduce that the first page of the spectral sequence is given by
\[ \Fpq{\sigma}{p}^{a,\bul}_1 \simeq H_c^{a,\bul+\dims\sigma+a}(\T^{\dims\sigma}) \otimes \SF^{p-a}(\ssinfty_\sigma). \]
Applying Proposition~\ref{prop:cohomology_compact_support_Tk}, we infer that $\Fpq{\sigma}{p}_1^{a,b}$ is trivial unless $a=\dims\sigma$ and $b=-\dims\sigma$, in which case it becomes equal to $\SF^{p-\dims\sigma}(\ssinfty_\sigma)$. Hence, $\Fpq{\sigma}{p}^{\bul,\bul}$ degenerates in page one, and $\Ep{p}^{\sigma,b}_1 \simeq \Tot^b(\Fpq{\sigma}{p}^{\bul,\bul}_1)$ is trivial unless $b=0$, in which case, we have $\Ep{p}^{\sigma,0}_1 \simeq \SF^{p-\dims\sigma}(\ssinfty_\sigma)$. This concludes the proof of the first part of Lemma~\ref{lem:hypercube_vanishing}. The second statement in the lemma follows from a cautious study of the different isomorphisms that we omit.  \qed

\section{Proof of Theorems~\ref{thm:ring_morphism} to~\ref{thm:Hodge_conjecture}} \label{sec:first_main_theorem}

We first prove Theorem~\ref{thm:ring_morphism}, and then deduce Theorems~\ref{thm:Hodge_isomorphism_matroids},~\ref{thm:Hodge_isomorphism},~\ref{thm:Hodge_isomorphism_dual} and~\ref{thm:Hodge_conjecture}. We assume that $\Sigma$ is unimodular, and saturated when required, and prove the theorems with integral coefficient. The proof of the statements in the theorems with rational coefficients will be similar.

\subsection{Vanishing part of Theorem~\ref{thm:ring_morphism}} \label{sec:vanishing}

By Theorem~\ref{thm:cohomology_compactification}, we have
\[H^{p,q}(\comp\Sigma) \simeq H^q(\cubC^{p,\bul}(\comp\Sigma)).\]
By the definition of the cubical complex \eqref{eq:definition_cubical_complex}, $\cubC^{p,q}(\comp\Sigma)$ is trivial for $p < q$, and therefore,
\[H^{p,q}(\comp\Sigma) =0 \qquad \textrm{ for } p<q.\]
For $q=0$ and $p>0$, we get
\[ H^{p,0}(\comp\Sigma) \simeq \ker\bigl(\SF^p(\conezero) \to \bigoplus_{\varrho \in \Sigma_1} \SF^{p-1}(\ssinfty_\varrho)\bigr). \]
Let $\alpha \in \SF^p(\conezero)$ be an element of the kernel. We show that $\alpha=0$. Let $V \subseteq \ssN_{\R}$ be the vector subspace spanned by $\Sigma$, that is, $V = \SF_{1}(\conezero)\otimes \R$. For any $\varrho \in \ssSigma_1$, we have $\contract_{\e_\varrho}(\alpha) = 0$. Let $\hat\alpha \in \bigwedge^p V^\dual$ be an element such that the restriction of $\hat\alpha$ to the subspace $\SF_p(\conezero) \subset \bigwedge^p V$ coincides with $\alpha$. Since the vectors $\e_\varrho$, $\varrho \in \Sigma_1$, span $V$, any element of $\bigwedge^p V^\dual$ whose contraction by all $\e_\varrho$ is trivial must be trivial. This means $\hat\alpha = 0$, which implies that $\alpha=0$, as required. We conclude that $H^{p,0}(\comp\Sigma) = 0$.

\subsection{Isomorphism between $H^{p,p}(\comp\Sigma)$ and $A^p(\Sigma)$}

In bidegree $(p,p)$, we get
\[ H^{p,p}(\comp\Sigma) \simeq \coker\bigl(\cubC^{p,p-1}(\comp\Sigma) \to \cubC^{p,p}(\comp\Sigma)\bigr) \simeq \coker\Bigl(\bigoplus_{\dims\tau=p-1}\SF^1(\ssinfty_\tau) \to \bigoplus_{\dims\sigma=p} \SF^0(\ssinfty_\sigma)\Bigr). \]
We have a natural isomorphism $\bigoplus_{\dims\sigma=p} \SF^0(\ssinfty_\sigma) \simeq Z^p(\Sigma)$, with $Z^p(\Sigma)$ the free abelian group generated by $\x_\sigma$ for $\sigma\in \Sigma_{p}$, defined in Section~\ref{sec:localization_lemma}. Moreover, since $\Sigma$ is saturated at any face $\tau\in\Sigma_{p-1}$, $\bigoplus_{\dims\tau=p-1}\SF^1(\ssinfty_\tau)$ generates the kernel of the surjective map $Z^p(\Sigma) \to A^p(\Sigma)$, as described in the Localization Lemma, Theorem~\ref{thm:localization-lemma}. Hence, we get an isomorphism $H^{p,p}(\comp\Sigma) \simeq A^p(\Sigma)$. By Lemma~\ref{lem:hypercube_vanishing}, the isomorphism is induced by the map
\[ \begin{tikzcd}[column sep=small, row sep=0em]
C^{p,p}(\comp\Sigma) \rar& \cubC^{p,p}(\comp\Sigma) \rar& Z^p(\Sigma) \\
a=(a_\delta)_{\substack{\delta\in\comp\Sigma \\ \dims\delta=p}} \ar[rr, mapsto] && \sum_{\sigma\in\Sigma_p}a_{\comp\sigma}(\ssnu_\sigma)\x_\sigma.
\end{tikzcd} \]

\subsection{Ring isomorphism}

It remains to prove that the isomorphism $H^{p,p}(\comp\Sigma) \simeq A^p(\Sigma)$ described in the previous section respects the products. This is more subtle and will be treated in Section~\ref{subsec:Ap_to_Hpp}. We will provide an explicit calculation of the inverse of the map $\Psi$ in Theorem~\ref{thm:ring_morphism} that shows that the cup-product in cohomology corresponds to product in the Chow ring.

\subsection{Proof of Theorem~\ref{thm:ring_morphism}}

Combining the results of the previous sections, we conclude. \qed

\subsection{Torsion-freeness and Proofs of Theorems~\ref{thm:Hodge_isomorphism_matroids} to~\ref{thm:Hodge_conjecture}} \label{sec:torsion_duality_chow}

Let $\Sigma$ be a unimodular fan. By the universal coefficient theorem applied to the tropical chain and cochain complexes of $\comp \Sigma$, for each pair of non-negative integers $p,q$, we get the following exact sequences
\begin{align}
\label{eq:univ_coeff_homology}
0\to \Ext{H^{p, q+1}(\comp\Sigma)}{\Z} \to H_{p,q}(\comp \Sigma) \to H^{p,q}(\comp \Sigma)^\dual \to 0,& \quad\text{and} \\
\label{eq:univ_coeff_cohomology}
0\to \Ext{H_{p, q-1}(\comp\Sigma)}{\Z} \to H^{p,q}(\comp \Sigma) \to H_{p,q}(\comp \Sigma)^\dual \to 0.&
\end{align}

From the first exact sequence, using the vanishing result stated in Theorem~\ref{thm:ring_morphism}, we obtain an isomorphism $H_{p,q}(\comp\Sigma) \simeq H^{p,q}(\comp\Sigma)^\dual$ for $p\leq q$. This shows the vanishing of $H_{p,q}(\comp \Sigma)$ for $p<q$ and an isomorphism $H_{p,p}(\comp \Sigma) \simeq H^{p,p}(\comp \Sigma)^\dual$. As a consequence, $H_{p,p}(\Sigma)$ is torsion-free.

\smallskip
\begin{proof}[Proof of Theorem~\ref{thm:Hodge_isomorphism_dual}] Let $\Sigma$ be a unimodular fan. The dual of the map from $H^{p,p}(\comp \Sigma)$ to $A^p(\Sigma)$ is the natural morphism from the group of Minkowski weights $\MW_p(\Sigma)$ (which is dual to $A^p(\Sigma)$ by Theorem~\ref{thm:duality_chow_mw}) to the tropical homology group $H_{p,p}(\comp\Sigma)$. When $\Sigma$ is saturated, Theorem~\ref{thm:ring_morphism} implies that this map is an isomorphism. By the discussion of Section~\ref{sec:ring_morphism_non_saturated}, this map is still an isomorphism even when $\Sigma$ is non-saturated. This combined with Proposition~\ref{prop:homology_q_0} stated below, proves Theorem~\ref{thm:Hodge_isomorphism_dual}.
\end{proof}

\begin{proof}[Proof of Theorem~\ref{thm:Hodge_conjecture}]
This is a direct consequence of Theorem~\ref{thm:Hodge_isomorphism_dual}.
\end{proof}

\begin{proof}[Proof of Theorem~\ref{thm:Hodge_isomorphism}]
In the case $\Sigma$ is a saturated unimodular tropical fan which is a tropical homology manifold with integral coefficients, we get an isomorphism $H^{p,q}(\comp \Sigma) \simeq H_{d-p,d-q}^\BM(\comp\Sigma) = H_{d-p,d-q}(\comp \Sigma)$, $p,q \in \Z_{\geq 0}$, which implies as well the vanishing of $H_{p,q}(\comp \Sigma)$ for $p>q$. Applying the exact sequence \eqref{eq:univ_coeff_cohomology}, we deduce the isomorphism $H^{p,q}(\comp \Sigma) \simeq H_{p,q}(\comp \Sigma)^\dual$. This in particular implies that the cohomology groups $H^{p,p}(\comp \Sigma)$, and so the Chow groups $A^p(\Sigma)$ in view of Theorem~\ref{thm:ring_morphism}, are all torsion-free. We infer that $A^p(\Sigma) \simeq \MW_p(\Sigma)^\dual$, and the Chow ring verifies Poincaré duality. Theorem~\ref{thm:Hodge_isomorphism} follows.
\end{proof}

\begin{proof}[Proof of Theorem~\ref{thm:Hodge_isomorphism_matroids}]
  For a matroid $\Ma$, the Bergman fan $\ssSigma_\Ma$ is unimodular, saturated, and a tropical homology manifold with integral coefficients by~\cite{JRS18}. Hence, applying Theorem~\ref{thm:Hodge_isomorphism}, we directly deduce that there is an isomorphism of rings $A^\bul(\Ma) \simto H^{\bul,\bul}(\comp\Sigma)$. Note that $\ssSigma_\Ma$ is a subfan of the permutahedral fan. Using this, it is easy to see that $\ss{\comp\Sigma}_\Ma$ is projective.
\end{proof}

It remains to prove the following proposition.

\begin{prop} \label{prop:homology_q_0}
The homology groups $H_{p,0}(\comp\Sigma)$ are trivial for each $p > 0$.
\end{prop}

A proof similar to the one given in Section~\ref{sec:vanishing} for the vanishing of $H^{p,0}(\comp\Sigma)$ shall give the result. Here, we propose a direct elementary argument which is valid without unimodularity and saturation property.

\begin{proof}
  Assume $p>0$. From the definition of the coefficient groups $\SF_p$, we see that for any $\tau\in\Sigma$, the map
  \begin{equation} \label{eqn:surj}
    \bigoplus_{\sigma \ssupface \tau} \SF_p(\cube_\sigma^\tau) \to \SF_p(\infty_\tau)
  \end{equation}
  is surjective. (If $\tau$ is a facet, $\SF_p(\infty_\tau)$ is trivial since $p>0$.) Take an element $a\in C_{p,0}(\comp\Sigma)$. By \eqref{eqn:surj} for $\tau=\conezero$, we find $b_0 \in C_{p,1}(\comp\Sigma)$ such that $a_1 \coloneqq a-\partial b_0$ is trivial at $\conezero$. Applying \eqref{eqn:surj} for $\tau=\rho$ with any ray $\rho \in \Sigma_1$, we get an element $b_1$ such that $a_2 \coloneqq a_1 - \partial b_1$ is trivial at $\conezero$ and also at $\infty_\rho$ for any ray $\rho\in\Sigma$. Proceeding by induction, we finally get that $a = \partial (b_0 + \dots + b_d)$ is a boundary. Hence, $H_{p,0}(\comp\Sigma)$ is trivial.
\end{proof}

\section{Explicit description of the inverse \texorpdfstring{$\Psi^{-1} \colon A^p(\Sigma) \to H^{p,p}(\comp \Sigma)$}{isomorphism}} \label{subsec:Ap_to_Hpp}

In this section, we provide an explicit description of the inverse of $\Psi$ which shows that the map $\Psi$ is a morphism of rings, concluding the proof of Theorem~\ref{thm:ring_morphism}. Following our convention, we assume that $\Sigma$ is unimodular and saturated, and work with integral coefficients. (Note however that the description we give in this section of the application $A^p(\Sigma) \to H^{p,p}(\comp \Sigma)$ works for any unimodular fan, see Section~\ref{sec:ring_morphism_non_saturated}.)

Before going through this, we introduce few notations that will help in working with the cohomology of $\comp\Sigma$. In the following, if $a$ is a cochain in $C^{p,q}(\comp\Sigma)$ and $\tau, \sigma \in \Sigma$ with $q=\dims\sigma-\dims\tau$, we denote by $\ssa^\tau_\sigma \in \SF^p(\cube^\tau_\sigma)$ the part of $a$ that lives on the face $\cube^\tau_\sigma$. If $\tau=\conezero$, we just write $a_\sigma$.

For any face $\delta=\cube^\tau_\sigma$ of dimension $q$ in $\comp\Sigma$, $\tau, \sigma\in \Sigma$ with $q=\dims\sigma-\dims\tau$, and for any $\alpha \in \SF^p(\delta)$, we denote by $[\delta, \alpha]$ the cochain in $C^{p,q}(\comp\Sigma)$ whose restriction to $\delta$ is $\alpha$ and which vanishes everywhere else.

\subsection{The inverse in cohomology degree one} \label{sec:inverse_degree_one}

We first consider the part $A^1(\Sigma)$. Let $\rho \in \Sigma_1$. We give an element of $H^{1,1}(\comp\Sigma)$ whose image by $\Psi$ coincides with $\x_\rho$. To this end, we first define an element $a \in C^{1,1}(\comp\Sigma)$ as follows. For any cone $\sigma\sim \rho$ in $\Sigma$, we choose an element $\ssalpha_\sigma \in \SF^1(\cube^{\sigma}_{\sigma \vee \rho})$ that takes value one on the vector $\sse^{\sigma}_{\sigma \vee \rho} \in \ssN^{\sigma}_{\infty, \sigma \vee \rho}$, and set
\[a \coloneqq \sum_{\sigma \sim \rho} [\cube^{\sigma}_{\sigma \vee \rho}, \ssalpha_\sigma].\]

Note that for any $\sigma \ssface \eta$, we have $\ssa^{\sigma}_\eta =0$ unless $\sigma \sim \rho$ and $\eta = \sigma \vee \rho$ in which case we have $\ssa^{\sigma}_\eta =\ssalpha_\sigma$. Moreover, the only part of $a$ that has sedentarity $\conezero$ is $[\cube_{\rho}, \ssalpha_\conezero]$, and we have $\Psi(a) =\ssalpha_\conezero(\sse_\rho)\x_\rho = \x_\rho$. The statement would follow if $a$ was a cocycle, which is not necessary the case in general. We will find an element $b$ in $C^{1,1}(\comp\Sigma)$ supported only on faces of non-zero sedentarity such that $a-b$ is a cocycle. We will then conclude by observing that since $\Psi(b)=0$, we still have $\Psi(a-b) =\x_\rho$, so that the inverse image of $\x_\rho$ will be represented by the class of $a-b$.

We describe $\hat a \coloneqq \d a$. Recall that $\sse^{\sigma}_{\sigma \vee \rho}$ is the projection at infinity of $\sse_\rho$ into $\ssN^\sigma_\infty$.
 Consider a pair of faces $\tau \subface \eta$ in $\Sigma$ with $\dims\eta-\dims\tau = 2$. We have
\[\sshata^\tau_\eta =0 \qquad \textrm{ if either } \rho \not\subface \eta \textrm{ or } \rho \subface \tau.\]
In the remaining cases, there must exist a ray $\rho'$ such that $\eta = \tau \vee \rho \vee \rho'$. In this case, we get
\[ \sshata^\tau_\eta = \pm(\i^*(\ssa^\tau_{\tau\vee\rho}) - \pi^*(\ssa^{\tau\vee\rho'}_{\tau\vee\rho'\vee\rho})) = \pm(\ssalpha_\tau - \pi^*(\ssalpha_{\tau\vee\rho'})), \]
where $\pi$ is the projection $\ssN^\tau \to \ssN^{\tau\vee\rho'}$.

Notice that $\sshata^\tau_\eta$ is zero on the projection of $\sse_\rho$. Therefore, we obtain a well-defined pushforward $\sspi^\rho_*(\sshata^\tau_\eta) \in \SF^1(\cube^{\tau\vee\rho}_\eta)$ where $\sspi^\rho\colon \ssN^\tau_\infty \to \ssN^{\tau\vee\rho}_\infty$ is the natural projection.

\smallskip
Set
\[ b \coloneqq \sum_{\tau,\eta \in \Sigma \\
\dims\eta-\dims\tau=2\\\rho \not \subface \tau, \, \rho \subface \eta} \sign(\cube^{\tau\vee\rho}_\eta, \cube^\tau_\eta) \Bigl[\cube^{\tau\vee\rho}_\eta,\sspi^{\rho}_{*}(\sshata^\tau_\eta)\Bigr]. \]
The element $b$ has been defined in order to get $\d b = \d a$ on every face not in $\compSigma^\rho_\infty$. In particular, the coboundary $c \coloneqq \d(a-b)$ of $a-b$ has support in $\compSigma^\rho_\infty$. Since $\d c=0$, the positive-sedentarity-vanishing Lemma~\ref{lem:psed_vanishing} below implies that $c =0$. This shows that $a-b$ is a cocycle.

Since $b$ has support only on the faces of non-zero sedentarity, we get in addition $\Psi(a-b) = \x_\rho$, as required.

\subsection{Positive-sedentarity-vanishing lemma}

We need the following lemma.
\begin{lemma}[Positive-sedentarity-vanishing lemma] \label{lem:psed_vanishing} Let $\Sigma$ be a unimodular fan. Let $\zeta$ be a non-zero face of\/ $\Sigma$ and let $c \in C^{p,q}(\comp \Sigma)$ be a cocycle supported in $\compSigma^\zeta_\infty$. Then we have $c=0$.
\end{lemma}

\begin{proof} By assumption, we have $\d c =0$. Since $c$ is supported in $\compSigma^\zeta_\infty$, we need to show the vanishing of $\ssc^\tau_\sigma$ for any pair $\tau, \sigma$ in $\Sigma$ with $\dims\sigma - \dims\tau =q$ and $\zeta\subface \tau$.

Let $\rho$ be a ray of $\zeta$, and consider the face $\eta \ssface \tau$ such that $\eta \wedge \rho = \conezero$ and $\tau = \rho \vee \eta$. Consider the face $\delta = \cube^\eta_\sigma$ of $\comp \Sigma$ and note that $\delta$ is of dimension $q+1$. The only face of dimension $q$ in $\delta$ that lies in $\compSigma^\zeta_\infty$ is $\cube^\tau_\sigma$. Since $c$ has support in $\compSigma^\zeta_\infty$, it follows that the components of $c$ on faces of $\delta$ different from $\cube^\tau_\sigma$ are all zero.

The projection map $\pi \colon \ssN^\eta_\infty \to \ssN^\tau_\infty$ induces a surjection $\sspi_*\colon \SF_p(\cube^\eta_\sigma) \to \SF_p(\cube^\tau_\sigma)$ and in injection $\sspi^*\colon \SF^p(\cube^\tau_\sigma)\to \SF^p(\cube^\eta_\sigma)$. We thus get
\[ \sspi^*(\ssc^\tau_\sigma) = \pm\ss{(\d c)}^\eta_\sigma =0.\]
We infer that $\ssc^\tau_\sigma =0$, and the lemma follows.
\end{proof}

\subsection{Cup-product in cubical complexes}

In order to extend the description of the inverse of $\Psi$ to higher degrees, we recall the formula for the cup product in cubical complexes. Let $a \in C^{p,q}(\comp\Sigma)$ and $b \in C^{p'\!\!,q'}(\comp\Sigma)$. The cup product $a \smile b$ is the element of $C^{p+p'\!\!,q+q'}(\comp\Sigma)$ with components described as follows. For any pair of faces $\tau \subface \eta$ with $\dims\eta - \dims\tau = q+q'$, the component $\ss{(a \smile b)}^\tau_\eta$ on the face $\cube^\tau_\eta$ is given by
\[ \ss{(a \smile b)}^\tau_\eta = \sum_{\tau\subface\sigma\subface\eta \\ \dims\sigma-\dims\tau = q} \ssvarpi^\tau_\eta(\ssnu^\tau_\sigma \wedge \sspi^*(\ssnu^\sigma_\eta))\ \cdot\ \i^*(\ssa^\tau_\sigma) \wedge \pi^*(\ssb^\sigma_\eta) \]
where, in the above sum, for the face $\tau \subface \sigma \subface \eta$, $\pi$ denotes the projection $\pi\colon \ssN^\tau \to \ssN^\sigma$.

As usual, this cup-product induces a cup product on cohomology
\[\smile\,\,\colon H^{p,q}(\comp\Sigma) \times H^{p',q'}(\comp\Sigma) \to H^{p+p',q+q'}(\comp\Sigma).\]

\subsection{The inverse in higher cohomological degrees}

{
\renewcommand{\a}[1]{{}_{\scaleto{#1}{4pt}}a}%
Let $\sigma \in \Sigma_p$ be a cone of dimension $p$ in $\Sigma$. By Localization Lemma, Theorem~\ref{thm:localization-lemma}, it suffices to find a preimage of $\x_\sigma$. Let $\ssrho_1, \dots, \ssrho_p$ be the rays of $\sigma$. For $i \in \{1,\dots,p\}$, we denote by $\a{i}$ the preimage of $\x_{\rho_i}$ as defined in Section~\ref{sec:inverse_degree_one}. Note in particular that $\a{i}$ is supported on the faces in $\comp\Sigma$ of the form $\cube^\tau_\sigma$ with $\ssrho_i\subface\sigma$.

We claim that the element
\[ a \coloneqq \a1 \smile \a2 \smile \dots \smile \a{p} \]
is a preimage of $\x_\sigma$.

We compute $a$ as follows. Denote by $\mathfrak S_p$ the symmetric group of order $p$. Consider a face $\sigma' \in \Sigma_p$ and denote by $\ssrho_1', \dots, \ssrho_p' \in \Sigma_1$ the rays of $\sigma'$. In the following, for a permutation $s \in \mathfrak S_p$ of $[p]$ and $k\in [p]$, we set
\[ \sigma'(s, k) \coloneqq \ssrho_{s(1)}' \vee \dots \vee \ssrho_{s(k)}', \]
the face of $\sigma'$ with rays $\ssrho_{s(1)}', \dots, \ssrho_{s(k)}'$.

\smallskip
Expanding the cup product using the formula stated in the previous section, we find
\[ \ssa_{\sigma'} = \sum_{s \in \mathfrak S_p} \ssvarpi_{\sigma'}(\ssnu_{\ssrho_{s(1)}'} \wedge \dots \wedge \ssnu_{\ssrho_{s(p)}'}) \ \cdot\ \ss{\a1}^{\sigma'(s,0)}_{\sigma'(s,1)} \wedge \ss{\a2}^{\sigma'(s,1)}_{\sigma'(s,2)} \wedge \dots \wedge \ss{\a{p}}^{\sigma'(s,p-1)}_{\sigma'(s,p)}, \]
where for the ease of reading, we omit to precise the pullback by different projections. Each term in the above sum is nontrivial only if $\rho_1 \subface \sigma'(s,1)$, $\rho_2 \subface \sigma'(s,2)$, \dots, and $\rho_p \subface \sigma'(s,p)$, \ie, if and only if $\sigma'=\sigma$ and the permutation $s$ is identity.

Since $\sspi_{\sigma(\id,j)}(\sse_{\rho_i}) = 0$ for $j \geq i$ and $\sspi_{\sigma(\id,j)}\colon N\to \ssN^{\sigma(\id,j)}$ the natural projection, we get that
\begin{align*}
\ssa_\sigma(\ssnu_\sigma)
  &= \ssvarpi_\sigma(\sse_{\rho_1}\wedge\dots\wedge\sse_{\rho_p}) \cdot \ssa_\sigma(\sse_{\rho_1}\wedge\dots\wedge\sse_{\rho_p}) \\
  &= \bigl(\ssvarpi_\sigma(\sse_{\rho_1}\wedge\dots\wedge\sse_{\rho_p})\bigr)^2 \bigl(\ss{\a1}^{\sigma(\id,0)}_{\sigma(\id,1)} \wedge \ss{\a2}^{\sigma(\id,1)}_{\sigma(\id,2)} \wedge \dots \wedge \ss{\a{p}}^{\sigma(\id,p-1)}_{\sigma(\id,p)}\bigr)\bigl(\sse_{\rho_1}\wedge\dots\wedge\sse_{\rho_p}\bigr) \\
  &= \ss{\a1}^{\sigma(\id,0)}_{\sigma(\id,0)\vee\rho_1}(\sse_{\rho_1})\ \cdots\ \ss{\a{p}}^{\sigma(\id,p-1)}_{\sigma(\id,p-1)\vee\rho_p}(\sse_{\rho_p}) \\
  &= 1.
\end{align*}
We thus infer that $a$ is a preimage of $\x_\sigma$. This achieves the description of the inverse of the map $\Psi$. By construction, the isomorphism respects the product.

\section{Characterization of tropical homology manifolds: Proof of Theorem~\ref{thm:smoothness_alternate}}\label{sec:proof-smoothness-criterion}

In Section~\ref{sec:fine_double_complex} we introduced the fine double complex $\Ep{p}^{\bul,\bul}$ and used the spectral sequence $\EpI{p}^{\bul,\bul}$ associated to the filtration by columns \eqref{sec:spectral_sequence_columns} to prove Theorem~\ref{thm:cohomology_compactification}. We consider now the spectral sequence $\EpII{p}^{\bul,\bul}$ given by the filtration by rows and use it to prove Theorem~\ref{thm:smoothness_alternate}.

By definition of the fine double complex, we know that the spectral sequence $\EpII{p}^{\bul,\bul}$ abuts to $H^{p,\bul}(\comp\Sigma)$
\[ \EpII{p}^{\bul,\bul} \Longrightarrow H^{p,\bul}(\comp\Sigma). \]

We will give the proof for unimodular tropical fans and the cohomology with integral coefficients. The proof for simplicial tropical fans and rational coefficients is similar.

We first prove the forward direction. Assume that the simplicial tropical fan $\Sigma$ is a tropical homology manifold with integral coefficients. Consider the $b$-th row of the double complex
\[ \EpII{p}_0^{\bul, b}\colon \quad \dots \longrightarrow
  \bigoplus_{\tau\subface\sigma \textrm{ faces of }\Sigma \\ \dims\sigma = a-1 \\ \dims\tau = -b} \SF^p(\cube^\tau_\sigma) \longrightarrow
  \bigoplus_{\tau\subface\sigma \textrm{ faces of }\Sigma \\ \dims\sigma = a \\ \dims\tau = -b} \SF^p(\cube^\tau_\sigma) \longrightarrow
  \bigoplus_{\tau\subface\sigma \textrm{ faces of }\Sigma \\ \dims\sigma = a+1 \\ \dims\tau = -b} \SF^p(\cube^\tau_\sigma) \longrightarrow
\cdots \]
Reorganizing $\EpII{p}_0^{\bul, b}$ as a sum according to $\tau$, we decompose $\EpII{p}_0^{\bul, b}$ as a direct sum of cochain complexes for the cohomology with compact support of the star fans of codimension $-b$:
\[ \EpII{p}_0^{\bul, b} = \bigoplus_{\dims\tau = -b} C_c^{p,\bul+b}(\Sigma^\tau). \]

Since $\Sigma$ is a tropical homology manifold with integral coefficients, all the star fans are tropical homology manifolds with integral coefficients, and so the cohomology with compact support $H_c^{p,q}(\Sigma^\tau,\Z)$ is trivial unless $q=d-(-b)$, and for this $q$, we have $H_c^{p,d+b}(\Sigma^\tau,\Z) \simeq H^{d+b-p,0}(\Sigma^\tau,\Z)^\dual = \SF_{d-p+b}(\ssinfty_\tau)$. Therefore, the first page of the spectral sequence is given by
\[ \EpII{p}_1^{a,b} = \begin{cases}
  \bigoplus_{\dims\tau = -b} \SF_{d-p+b}(\ssinfty_\tau) & \text{if $a = d$,} \\
  0 & \text{otherwise.}
\end{cases} \]
In other words, $\EpII{p}_1^{a,b}$ is trivial unless $a=d$ and we have
\[ \EpII{p}_1^{d,\bul} = \hcubC_{d-p,-\bul}(\comp\Sigma) \]
where $\hcubC_{d-p,\bul}(\comp\Sigma)$ denotes the dual to the complex $\cubC^{d-p,\bul}(\comp\Sigma)$ introduced in~\eqref{eq:definition_cubical_complex}. Theorem~\ref{thm:cohomology_compactification} combined with the universal coefficient theorem imply that the homology of $\hcubC_{d-p,\bul}(\comp\Sigma)$ computes $H_{d-p,\bul}(\comp\Sigma)$. Hence, $\EpII{p}_1^{\bul,\bul}$ degenerates at page two, and we have
\[ \EpII{p}^{\bul,\bul} \Longrightarrow H_{d-p,d-\bul}(\comp\Sigma). \]
Since
\[ \EpII{p}^{\bul,\bul} \Longrightarrow H^{p,\bul}(\comp\Sigma). \]
we get the isomorphism of Poincaré duality $H^{p,q}(\comp\Sigma,\Z) \simeq H_{d-p,d-q}(\comp\Sigma,\Z)$. More precisely, we infer that
\[ H^{p,p}(\comp\Sigma) \simeq \EpI{p}_\infty^{p,0} \simeq \EpII{p}_\infty^{d,p-d} \simeq H_{d-p,d-p}(\comp\Sigma), \]
and all the other terms in $\EpI{p}^{\bul,\bul}_\infty$ and $\EpII{p}^{\bul,\bul}_\infty$ are trivial.

Applying the same argument to all the star fans $\Sigma^\sigma$, $\sigma\in\Sigma$, we conclude with the proof of the forward direction.

\smallskip
We now prove the other direction. Assume that $\comp\Sigma^\sigma$ verifies Poincaré duality for any face $\sigma$ of $\Sigma$. By induction, we can assume that $\Sigma^\sigma$ is a tropical homology manifold with integral coefficients for any nonzero face $\sigma$ of $\Sigma$.

We have a natural inclusion of complexes $C_c^{p,\bul}(\Sigma) \simeq \Ep{p}^{\bul,0} \hookrightarrow \Ep{p}^{\bul,\bul}$. Denote by $\Ep{p}^{\bul,\bul<0}$ the cokernel. We get a short exact sequence of complexes:
\begin{equation} \label{eqn:cut_E}
0 \longrightarrow C_c^{p,\bul}(\Sigma) \longrightarrow \Tot^{\bul}(\Ep{p}^{\bul,\bul}) \longrightarrow \Tot^{\bul}(\Ep{p}^{\bul,\bul<0}) \longrightarrow 0.
\end{equation}

We already described the cohomology of the first two terms in the sequence. For the third one, we compute the first page of the spectral sequence associated to filtration by columns, use the induction hypothesis that the proper star fans are tropical homology manifolds, and get
\[\EpI{p}_1^{d,\bul<0} = \hcubC_{d-p,-\bul>0}(\comp\Sigma)\]
where
\[ \hcubC_{d-p,\bul>0}(\comp\Sigma) \coloneqq \coker\Bigl(\SF_{d-p}(\conezero) \hookrightarrow \hcubC_{d-p,\bul}(\comp\Sigma)\Bigr). \]
Moreover, all the other terms in the spectral sequence $\EpI{p}_1^{\bul,\bul<0}$ are trivial.

We treat first the case $p < d-1$. The cases $p=d-1$ and $p=d$ will be studied later below.

In this case, we get
\[ H_k(\hcubC_{d-p,\bul>0}(\comp\Sigma)) = \begin{cases}
  H_{d-p,d-p}(\comp\Sigma) & \text{if $k = d-p$,} \\
  \SF_{d-p}(\conezero)     & \text{if $k = 1$,} \\
  0                        & \text{otherwise.}
\end{cases} \]

The long exact sequence associated to \eqref{eqn:cut_E} leads to
\[ \cdots \longrightarrow H_c^{p,q}(\Sigma) \longrightarrow H^{p,q}(\comp\Sigma) \longrightarrow H^q(\hcubC_{d-p,(d-\bul)>0}(\comp\Sigma)) \longrightarrow H_c^{p,q+1}(\Sigma) \longrightarrow \cdots \]
The two middle terms are trivial unless $q = p$ or $q = d-1$.

In the case $q=p$, we get a map $H^{p,p}(\comp\Sigma) \to H_{d-p,d-p}(\comp\Sigma)$. By assumption, this is an isomorphism.

For $q=d-1$, since $p<d-1$, we have $H^{p,d-1}(\comp\Sigma) = H^{p,d}(\comp\Sigma) = 0$, and we get
\[ 0 \longrightarrow \SF_{d-p}(\conezero) \longrightarrow H_c^{p,d}(\Sigma) \longrightarrow 0, \]
that is, an isomorphism $\SF_{d-p}(\conezero) \simeq H_c^{p,d}(\Sigma)$.

All the other terms of $H_c^{p,\bul}(\Sigma)$ vanish.

Using the universal coefficient theorem, we get $H^\BM_{p,d}(\Sigma) \simeq H_c^{p,d}(\Sigma)^\dual \simeq \SF^{d-p}(\conezero)$, and $H^\BM_{p,q}(\Sigma) = 0$ for $q \neq d$. This is Poincaré duality for the homology with coefficients in~$\SF_p$, see~\ref{subsec:homology}.

\smallskip
In the case $p=d$, the cohomology $H_\bul(\hcubC_{0,\bul>0}(\comp\Sigma))$ is trivial. Hence we have an isomorphism $H_c^{d,\bul}(\Sigma) \simeq H^{d,\bul}(\comp\Sigma)$ and we conclude in the same way as above.

\smallskip
It remains to treat the case $p = d-1$. In this case, $\hcubC_{1,\bul>0}(\Sigma)$ contains only one nontrivial term: $\hcubC_{1,1>0}(\Sigma) = \bigoplus_{\varrho\in\Sigma_1} \SF_0(\ssinfty_\varrho)$. Its cohomology is identical. In the above long exact sequence, the only interesting part is the first row of the following diagram
\[ \begin{tikzcd}
0 \rar& H_c^{d-1,d-1}(\Sigma) \rar\dar& H^{d-1,d-1}(\comp\Sigma) \rar\dar{\vsim}& \bigoplus_{\varrho\in\Sigma_1}\SF_0(\ssinfty_\varrho) \rar\dar{\vsim}& H_c^{d-1,d}(\Sigma) \rar\dar& 0 \\
 & 0 \rar& \MW_1(\Sigma) \rar& \W_1(\Sigma) \rar& \SF_1(\conezero) \rar& 0
\end{tikzcd} \]
We describe the second row and the vertical maps. We know that $H^{d-1,d-1}(\comp\Sigma) \simeq H_{1,1}(\comp\Sigma) \simeq \MW_1(\Sigma)$. Here, $\W_1(\Sigma) \coloneqq \Z^{\Sigma_1} \simeq \bigoplus_{\varrho\in\Sigma_1}\SF_0(\ssinfty_\varrho)$ denotes the group of one dimensional weights of $\Sigma$. Note also that we have a natural map $H_c^{d-1,d}(\Sigma) \to \SF_1(\conezero)$ described in Section~\ref{sec:tropical-fans}. By the definition of the Minkowski weights, the second row is a short exact sequence. The five lemma implies that the first and the last vertical maps are isomorphisms. Therefore, $H_c^{d-1,q}(\Sigma)$ is trivial for $q\neq d$, and $H_c^{d-1,q}(\Sigma) \simeq \SF_1(\conezero)$. Using the universal coefficient theorem as above, and dualizing, we get Poincaré duality for $H^\BM_{d-1,\bul}(\Sigma)$. Altogether we obtain Poincaré duality for $H^\BM_{\bul,\bul}(\Sigma)$. \qed

\section{Proof of Theorem~\ref{thm:positivity}} \label{sec:kleiman}

In this section, we prove Theorem~\ref{thm:positivity}. Let $\alpha$ be a class in $A^1(\Sigma)$ associated to a conewise linear function $f$ on $\Sigma$, that is,
\[\alpha = \sum_{\rho\in \Sigma_1} f(\sse_\rho)\ssx_\rho \in A^1(\Sigma)\]
where $\ssx_\rho$ is the class of $\x_\rho$ in $A^1(\Sigma)$. Via Theorem~\ref{thm:ring_morphism}, we identify $\alpha$ as an element of $H^{1,1}(\Sigma)$.

We first prove the forward direction. Assume that $\alpha$ is ample in $A^1(\Sigma)$. This means that the function $f$ is strictly convex on $\Sigma$. Then, $f$ induces a strictly convex function $\ssf^\sigma$ on each stratum $\ssSigma_\infty^\sigma$ (well-defined up to a linear function). We can assume that $\ssf^\sigma$ is positive away from the origin $\ssinfty_\sigma$ of the fan $\ssSigma_\infty^\sigma$. On any ray $\rho=\sseta^\sigma_\infty$ of $\comp\Sigma$, we have $\ssf^\sigma(\e^\sigma_\eta) > 0$. For any effective nonzero Minkowski weight $\gamma = (\ssSigma_{\infty,1}^\sigma, w)$ of dimension one supported on $\ssSigma_\infty^\sigma$, we thus get
\[\langle\alpha,\gamma \rangle = \sum_{\eta\ssupface \sigma} w(\sseta^\sigma_\infty)\ssf^\sigma(\e_\eta^\sigma)>0\]
proving the forward direction for effective tropical curves in $\comp\Sigma$ which are given by an effective nonzero Minkowski weight on some stratum of $\comp\Sigma$. The proof can be adapted to give positivity of the pairing with more general effective nonzero tropical curves in $\comp\Sigma$. We omit the details.

\begin{remark}
  Note that it follows from the theorem that the effective cone in $H^{1,1}(\comp\Sigma)$ is generated by effective Minkowsi weights of dimension one in strata. A similar statement can be proved in any dimension.
\end{remark}

For the other direction, proceeding by induction, we can assume that the restriction $\ssalpha^\sigma$ of $\alpha$ to $\ssSigma^\sigma_\infty$ is ample for any $\sigma \neq \conezero$. Therefore, $f$ is strictly convex at any cone $\sigma \neq \conezero$ of $\Sigma$. We show that it is also strictly convex at $\conezero$.

The following argument is due to Pierre-Louis Blayac. Consider the restriction $f\rest{\supp{\ssSigma_{(1)}}}$ of $f$ to the one-skeleton of $\Sigma$. Take the graph $\Gamma$ of $f\rest{\supp{\ssSigma_{(1)}}}$ in $\ssN \times \R$. Let $C$ be the convex hull of $\Gamma\,\cup\,(\{0\}\times \R_+)$. The origin is in $C$. Two cases can occur.
\begin{itemize}
  \item Either, $C$ contains a point of $\{0\}\times\R_{<0}$. We find a collection $(c_\rho)_{\rho\in \Sigma_1}$ of nonnegative numbers such that
  \[\sum_{\rho \in \Sigma_1} c_\rho(\e_\rho, f(\e_\rho)) = (0, -1).\]
  We have found a Minkowski weight $\gamma \in \MW_1(\Sigma, \R)$ given by the reals $c_\rho$ such that $\langle \alpha, \gamma\rangle <0$. Since $\MW_1(\Sigma, \Q)$ is dense in $\MW_1(\Sigma, \R)$, this implies the existence of an element $\theta\in \MW_1(\Sigma)$ with $\langle\alpha,\theta\rangle <0$, leading to a contradiction.
  \item Or, there is a support hyperplane $H$ of $C$ at $0$ separating $\{0\}\times\R_{<0}$ from $C$. If $H\cap C$ is strictly convex, we can move $H$ a bit such that $H\cap C=\{0\}$. In this case, we conclude that $f$ is strictly convex at $\conezero$, as required. Otherwise, $H\cap C$ contains a non-zero vector subspace $V$. The space $V$ is positively generated as the sum of rays of $\Gamma$ that lie in $V$. We thus get an effective element $\gamma$ of $\MW_1(\Sigma)$ with $\langle\alpha, \gamma\rangle=0$, leading again to a contradiction.
\end{itemize}

\section{Examples and discussions} \label{sec:discussion}

In this section, we provide examples and discuss complementary results.

\subsection{Torsion in the Chow ring of non-saturated fans and Theorem~\ref{thm:ring_morphism}} \label{sec:ring_morphism_non_saturated}

The proof of Theorem~\ref{thm:ring_morphism} gives results in the case where $\Sigma$ is unimodular but not saturated. In this situation, we still have a surjective ring morphism $A^\bul(\Sigma) \to H^{\bul,\bul}(\comp\Sigma)$ with a kernel that is torsion and can be nontrivial, see Example~\ref{ex:F1_vs_N} below. On the other hand, the kernel can be trivial even though $\Sigma$ is not saturated, see \cite[Example 12.14]{AP-hodge-fan}. In any case, as stated in Theorem~\ref{thm:ring_morphism}, we have an isomorphism for rational coefficients, as well as a dual isomorphism with integral coefficients: $A^\bul(\Sigma)^\dual \simeq H^{\bul,\bul}(\comp\Sigma)^\dual$. The part of the theorem that concerns the vanishing of cohomology groups for $p<q$ is still valid for any unimodular fan without the saturation hypothesis, and for $p>q=0$ it is valid in full generality.

\begin{example} \label{ex:F1_vs_N}
  We recall the example \cite[Example 12.11]{AP-hodge-fan}. Let $(\e_1,\e_2)$ be the standard basis of $\Z^2$ and let $\e_0=-\e_1-\e_2$. Denote by $\rho_i = \R_{\geq 0}\e_i$, $i\in\{0,1,2\}$ the corresponding rays.

  \smallskip
  In the following, we consider the lattice $N \coloneqq \Z\e_1 + \frac 13\Z(\e_1 - \e_2)$ in $\R^2$. Note that $\Z^2 \subset N$ is a sublattice of index three in $N$. The dual lattice $M \coloneqq N^\dual$ is of index three in $(\Z^2)^\dual$.

  Let $\Delta$ be the one-dimensional fan in $N$ with rays $\rho_0, \rho_1, \rho_2$ and lattice $\ssN_\Delta = N$. The fan $\Delta$ is tropical and unimodular but it is not saturated. Any element $f$ of $(\Z^2)^\dual \setminus M$ induces a meromorphic function on $\Delta$ which is linear but not integral linear. (The meromorphic function $3f$ will be integral linear.) As a consequence, $A^1(\Delta)$ has torsion: the element $x_{\rho_1}-x_{\rho_2}$ is non-zero, but $3\cdot(x_{\rho_1}-x_{\rho_2})$ vanishes in $A^1(\Delta)$.

  However, $\SF_1(\conezero)$ is by definition equal to $\Z^2\subseteq N$. Hence, $\SF^1(\conezero) = (\Z^2)^\dual$, and a direct computation shows that $H^{1,1}(\comp\Delta) \simeq \Z$. The map from $A^1 (\Delta)$ to $H^{1,1}(\comp\Delta)$ is not an isomorphism. It is surjective with kernel equal to the torsion part of $A^1(\Delta)$.
\end{example}

\subsection{Torsion in integral cohomology of non-unimodular fans}\label{sec:non_unimodular_projective_fan}

The following example shows that the unimodularity assumption is in general needed in the second part of Theorem~\ref{thm:ring_morphism}.

\begin{example} \label{ex:non_unimodular_projective_fan}
  Let $\Sigma$ be the complete fan in the plane $\R^2$, with $N =\Z^2$, whose rays are the $\rho_i$ defined as in Example~\ref{ex:F1_vs_N}. The fan $\Sigma$ is not unimodular. As for Example~\ref{ex:F1_vs_N}, $A^1(\Sigma)$ has torsion. Moreover, the cohomology group $H^{1,2}(\comp\Sigma) \simeq \rquot\Z{3\Z}$ is nontrivial even though we are in bidegree $(p,q) = (1,2)$ with $p < q$.
\end{example}

\subsection{Necessity of being tropical homology manifold in Theorem~\ref{thm:Hodge_isomorphism}}

The assumption made on $\Sigma$ being a tropical homology manifold in the statement of Theorem~\ref{thm:Hodge_isomorphism} is needed, and the result is not true in general, even for a saturated unimodular tropical fan. Example~\ref{ex:cube} is a saturated unimodular tropical fan $\Sigma$ with $H^{2,1}(\comp\Sigma)$ of rank two.

\begin{example}[The fan over the one-skeleton of the cube] \label{ex:cube}
We consider the fan defined over the $1$-skeleton of a cube, thoroughly discussed in our work~\cite[Section 12.3]{AP-hodge-fan}. We consider the standard cube $\mbox{\mancube}$ with vertices $(\pm1, \pm1, \pm1)$, and let $\Sigma$ be the two-dimensional fan with rays generated by vertices and with facets generated by edges of the cube. The fan $\Sigma$ is locally irreducible and tropical but not unimodular. We obtain a unimodular fan by changing the underlying lattice. In what follows, we work with the lattice $N := \sum_{\varrho\in\Sigma_1}\Z\e_\varrho$.

Table~\ref{tab:cohomology_cube} summarizes the main cohomological data about the cube. We can compute the homology via the universal coefficient theorem: the torsion part is shifted by one column on the left, the rest remains unchanged. Moreover the image of $A^1(\Sigma)$ inside $A^1(\Sigma)^\dual$ is a sublattice of full rank of index two.
\end{example}

\begin{table}
\renewcommand{\strut}{\rule{0pt}{1.1em}}
\[ \begin{array}{c|c|c|c|cc|c|c|c|}
\cline{2-4} \cline{7-9}
\strut H^{p,q}_c(\Sigma)       & q = 0 & 1    & 2                         & \qquad\qquad          & H^{p,q}(\comp\Sigma) & q = 0  & 1    & 2    \\ \cline{1-4} \cline{6-9}
\multicolumn{1}{|c|}{\strut p = 0} & 0 & 0    & \Z^5                      & \multicolumn{1}{c|}{} & p = 0                    & \Z & 0    & 0    \\ \cline{1-4} \cline{6-9}
\multicolumn{1}{|c|}{\strut 1} & 0 & 0    & \Z^3\times\rquot{\Z}{2\Z} & \multicolumn{1}{c|}{} & 1                    & 0  & \Z^5 & 0    \\ \cline{1-4} \cline{6-9}
\multicolumn{1}{|c|}{\strut 2} & 0 & \Z^2 & \Z                        & \multicolumn{1}{c|}{} & 2                    & 0  & \Z^2 & \Z   \\ \cline{1-4} \cline{6-9}
\end{array} \]
\caption{Cohomology of the fan $\Sigma$ over the one-skeleton of the cube. \label{tab:cohomology_cube}}
\end{table}

\subsection{Explicit description of Poincaré duality}

Let $\Sigma$ be a unimodular fan of dimension $d$. Assume that $\Sigma$ is a tropical homology manifold with integral coefficients. By Poincaré duality, there is an isomorphism $\PD\colon H^{p,p}(\comp\Sigma) \simto H_{d-p,d-p}(\comp\Sigma)$. We describe this map in two different ways, as follows.

Let $\alpha \in H^{p,p}(\comp\Sigma)$ and choose a representative $a = (a_\delta)_{\delta\in\comp\Sigma, \dims\delta=p}$ of $\alpha$ in $C^{p,p}(\comp\Sigma)$. For any $\sigma \in \Sigma$, denote by $[\compSigma_\infty^\sigma] \in H_{d-\dims\sigma,d-\dims\sigma}(\comp\Sigma)$ the image of the fundamental class of $\compSigma^\sigma_\infty$ in the homology of $\comp\Sigma$. Then, we have
\[ \PD(\alpha) = \sum_{\sigma\in \Sigma_p} a_{\comp\sigma}(\ssnu_\sigma) [\compSigma^\sigma_\infty]. \]

Alternatively, we describe $\PD(\alpha)$ as a Minkowski weight, through the cycle class map $\MW_{d-p}(\Sigma) \to H_{d-p,d-p}(\comp\Sigma)$. Let $w\colon \Sigma_{d-p} \to \Z$ be the map defined by $w(\sigma) \coloneqq \langle \alpha, [\compSigma^\sigma_\infty] \rangle$. Then, $w$ is a Minkowski weight and the corresponding cycle in $H_{d-p,d-p}(\comp\Sigma)$ is $\PD(\alpha)$.

\subsection{The case of non-rational simplicial fans} \label{sec:nonrational}

The results of this paper remain valid for simplicial fans which are not necessarily rational, working with real coefficients. Let $\Sigma$ be a simplicial fan in $\ssN_{\R}$. For each ray $\rho \in \Sigma_1$, we choose a generator $\sse_\rho$ so that $\R_+\sse_\rho = \rho$. We define the Chow ring $A^\bul(\Sigma, \R)$ of $\Sigma$ as quotient
\[ A^\bul(\Sigma, \R) \coloneqq \rquot{\R[\x_\zeta]_{\zeta\in \Sigma_1}}{\bigl(I + J\bigr)} \]
where $I$ is the ideal generated by the products $\x_{\rho_1}\!\cdots \x_{\rho_k}$, for $k\in \N$, such that $\ssrho_1, \dots, \ssrho_k$ are non-comparable rays in $\Sigma$, and $J$ is the ideal generated by the elements of the form
\[ \sum_{\zeta\in \Sigma_1} m(\sse_\zeta)\x_\zeta, \qquad  m \in \ssM_\R.\]
Different choices of the generators $\sse_\rho$ for $\rho\in \Sigma_1$ give isomorphic Chow rings.

Similarly, we define the sheaves $\SF^p(\,\cdot\,, \R)$ with real coefficients, and the cohomology groups $H^{p,q}(\comp \Sigma, \R) \coloneqq H^q(\comp \Sigma, \SF^p(\,\cdot\,, \R))$.

All the theorems stated in the introduction remain valid with the above adjustment of the coefficients.

\printbibliography

\end{document}